\documentclass[10pt]{article}
\usepackage{amsmath,amssymb,amsfonts,mathrsfs,latexsym,parskip,mathtools}
\usepackage{graphicx, subfigure}
\graphicspath{{Figures/}}
\usepackage[justification=justified,singlelinecheck=false]{caption}
\usepackage[hidelinks]{hyperref} % Hyperlinks for equation numbers
\usepackage[text={17cm,23cm},top=3cm]{geometry}  % This line determines the page layout.
\usepackage{titlesec}
\usepackage[title]{appendix}

\numberwithin{equation}{section}

\pdfoutput=1

\usepackage[usenames,dvipsnames,svgnames]{xcolor}

\hypersetup{
    colorlinks,
    linkcolor={blue!50!black},
    citecolor={blue!50!black},
    urlcolor={blue!50!black}
}

\usepackage{tabularx}
\usepackage{hhline}
\newcolumntype{Y}{>{\centering\arraybackslash}X}

\usepackage{nth}

\usepackage{color}

\newcolumntype{b}{X}
\newcolumntype{s}{>{\hsize=.5\hsize}X}

\usepackage{afterpage}
\makeatletter
\newcommand{\oneraggedpage}{\let\mytextbottom\@textbottom
  \let\mytexttop\@texttop
  \raggedbottom
  \afterpage{%
  \global\let\@textbottom\mytextbottom
  \global\let\@texttop\mytexttop}}
  
\allowdisplaybreaks[1]

\usepackage{comment}

\newcommand{\lb}{\left (}
\newcommand{\rb}{\right )}
\newcommand{\lset}{\left \{}
\newcommand{\rset}{\right \}}
\newcommand{\eqtext}[1]{\quad \text{#1} \quad}
\newcommand{\lsq}{\left [}
\newcommand{\rsq}{\right ]}

\newcommand{\sechn}[2]{\mathrm{sech}^{#1} \lb #2 \rb}
\newcommand{\tanhn}[2]{\mathrm{tanh}^{#1} \lb #2 \rb}
\newcommand{\dd}[1]{\; \mathrm{d} #1}
\newcommand{\diff}[2]{\frac{\mathrm{d} #1}{\mathrm{d} #2}}
\newcommand{\diffn}[3]{\dfrac{\mathrm{d}^{#1} #2}{\mathrm{d} #3^{#1}}}

\newcommand{\pinv}{\partial^{-1}}

\renewcommand{\O}[1]{\mathcal{O} \lb #1 \rb}

\begin{document}

% Symbolic footnotes
\renewcommand*{\thefootnote}{\fnsymbol{footnote}}

% Front matter
\begin{center}
{\LARGE Periodic solutions of coupled Boussinesq equations and Ostrovsky-type models free from zero-mass contradiction} \\[2.5em]
{\large K. R. Khusnutdinova$^{a}$, M. R. Tranter$^{b}$ \footnote{Corresponding author. Tel: +44 (0)1158 483412.}} \\[1.5em]
{$^{a}$ Department of Mathematical Sciences, Loughborough University, Loughborough LE11 3TU, UK} \\[1.5em]
{$^{b}$ Department of Physics and Mathematics, Nottingham Trent University, Nottingham NG11 8NS, UK} \\[1.5em]
{K.Khusnutdinova@lboro.ac.uk} \\
{Matt.Tranter@ntu.ac.uk}
\end{center}

% Arabic footnotes
\renewcommand*{\thefootnote}{\arabic{footnote}}
\setcounter{footnote}{0}

% Abstract
\begin{abstract}
Coupled Boussinesq equations describe long weakly-nonlinear longitudinal strain waves in a bi-layer with a soft bonding between the layers (e.g. a soft adhesive). From the mathematical viewpoint, a particularly difficult case appears when the linear long-wave speeds in the layers are significantly different (high-contrast case). The traditional derivation of the uni-directional models leads to four uncoupled Ostrovsky equations, for the right- and left-propagating waves in each layer.  However, the models impose a ``zero-mass constraint'' i.e. the initial conditions should necessarily have zero mean, restricting the applicability of that description. Here, we bypass the contradiction in this high-contrast case by constructing the solution for the deviation from the evolving mean value, using asymptotic multiple-scale expansions involving two pairs of fast characteristic variables and two slow-time variables. By construction, the Ostrovsky equations emerging within the scope of this derivation are solved for initial conditions with zero mean while initial conditions for the original system may have non-zero mean values. Asymptotic validity of the solution is carefully examined numerically. We apply the models to the description of counter-propagating waves generated by solitary wave initial conditions, or co-propagating waves generated by cnoidal wave initial conditions, as well as the resulting wave interactions, and contrast with the behaviour of the waves in bi-layers when the linear long-wave speeds in the layers are close (low-contrast case). 
One local (classical) and two non-local (generalised) conservation laws of the coupled Boussinesq equations for strains are derived, and these are used to control the accuracy of the numerical simulations.
\end{abstract}

%\maketitle

\section{Introduction}
Korteweg-de Vries and Boussinesq-type equations have been derived to describe long weakly-nonlinear longitudinal strain waves in rod- and bar-like elastic solids  \cite{Nariboli70, Ostrovsky77, Samsonov01, Porubov03, Erofeev02, Dai04, Peets17, Garbuzov19, Garbuzov20}. This paper will focus on the system of coupled regularised Boussinesq (cRB) equations, presented here in non-dimensional and scaled form:
\begin{align}
	u_{tt} - u_{xx} &= \varepsilon \lsq \frac{1}{2} \lb u^2 \rb_{xx} + u_{ttxx} - \delta \lb u - w \rb \rsq, \label{ueq} \\
	w_{tt} - c^2 w_{xx} &= \varepsilon \lsq \frac{\alpha}{2} \lb w^2 \rb_{xx} + \beta w_{ttxx} + \gamma \lb u - w \rb \rsq. \label{weq}
\end{align}
This system of equations describes long nonlinear longitudinal strain waves in a bi-layer with a soft bonding between the layers, allowing the layers to move relative to each other \cite{Khusnutdinova09}. In this context, $u$ and $w$ denote longitudinal strains in the layers, $\alpha$, $\beta$, $\gamma$, $\delta$ are coefficients depending on the mechanical and geometrical properties of a waveguide, $c$ is the ratio of the characteristic linear wave speeds in the layers and $\varepsilon$ is a small amplitude parameter. When the layers have similar properties, radiating solitary waves are found to propagate on large periodic domains \cite{Khusnutdinova17}.

In all of these contexts, the natural initial conditions have non-zero mean, and the associated uni-directional equations emerging in the construction of weakly-nonlinear solutions are typically solved numerically using pseudo-spectral schemes with periodic boundary conditions. However, recent work has shown that solving such equations on a periodic domain may require careful attention if the initial conditions have non-zero mean and the periodic domain is comparable to the scale of the initial condition \cite{Khusnutdinova14, Khusnutdinova19}. This was investigated for the Boussinesq-Klein-Gordon (BKG) equation, which can be obtained from the cRB equations by taking the limit $\gamma \rightarrow 0$ with $w = 0$. The weakly-nonlinear solution derived via the traditional procedure leads, at leading order, to two uni-directional Ostrovsky equations, originally developed in the context of fluids \cite{Ostrovsky78} (see also \cite{Grimshaw98}). The model necessarily requires that regular solutions have zero mean initial conditions. The existence of such a formal constraint is known as the ``zero-mass (or zero mean) contradiction'' since the original problem formulation does not impose this constraint \cite{Khusnutdinova19}. However, it was shown that the derivation procedure can be modified in order to develop a weakly-nonlinear solution of the BKG equation for a deviation from the evolving non-zero mean \cite{Khusnutdinova19}. Earlier results were developed at the level of Fourier expansions in the spatial variable \cite{Khusnutdinova14}, while the procedure suggested in \cite{Khusnutdinova19} can be viewed as a nonlinear extension of the d'Alembert solution. The Ostrovsky equations emerging within the scope of this procedure are solved for zero mean initial conditions by construction, and the zero mean contradiction is avoided. We note that additional conservation laws of the ``moment'' type were constructed for the Ostrovsky equation in \cite{Benilov92} (see also \cite{Stepanyants20}). However, these conservation laws are not applicable in our case since multiplying a periodic function of some variable by powers of that variable takes us outside of the class of periodic functions. Hence, such conservation laws impose no additional constraints in the periodic case under study.

Considering our cRB equations, when the linear characteristic speeds are close, satisfying the relation $c - 1 = \mathcal{O}(\varepsilon)$ (low-contrast case), and the period of the solution is large compared to the scale of a localised initial condition, we find long-living radiating solitary waves \cite{Khusnutdinova09, Khusnutdinova11, Grimshaw17}. When $c - 1 = \mathcal{O}(1)$, we find wave packets governed by the Ostrovsky equations \cite{Khusnutdinova11, Ostrovsky78}. In this paper we consider the case when $c - 1 = \mathcal{O}(1)$ (high-contrast case), and the period of the solution is comparable with the scale of the localised initial condition, or when the initial condition is not localised at all.  Previously, a simpler case when $c - 1 = \mathcal{O}(\varepsilon)$ was investigated in \cite{Khusnutdinova19Book}. As we will obtain single or coupled Ostrovsky equations to leading order, depending on the assumption on the characteristic speeds in the layers,  we will need to consider how a weakly-nonlinear solution can be constructed that takes account of the zero-mass contradiction. 

%%%%%%
The paper is organised as follows. In Section \ref{sec:WNL} we construct a weakly-nonlinear solution of the Cauchy problem for the cRB equations (\ref{ueq}) - (\ref{weq}) in the case $c - 1 = \mathcal{O}(1)$ on a periodic domain, using asymptotic multiple-scale expansions for the deviation from the oscillating mean values. We use two sets of fast characteristic variables and two slow time variables. The validity of the solutions is examined in Section \ref{sec:Validity} by comparing the constructed weakly-nonlinear solution with direct numerical simulations. 
In Section \ref{sec:Counter} we use both direct numerical simulations and the constructed weakly-nonlinear solutions to study the interaction of counter-propagating waves generated by solitary wave initial conditions, or co-propagating waves generated by cnoidal wave initial conditions, when the characteristic speeds are either close or significantly different. We only consider co-propagating waves for cnoidal wave initial conditions, as this case provides a venue for the study of the strong wave interactions. We also determine the nature of the interaction. 
In Section \ref{sec:Laws} we discuss the local and non-local conservation laws used to control the accuracy of numerical simulations and we conclude our studies in Section \ref{sec:Conc}. The numerical schemes used in our numerical simulations are contained within the appendices.

\section{Weakly Nonlinear Solution}
\label{sec:WNL}
We solve the equation system \eqref{ueq} - \eqref{weq} on the periodic domain $x \in [-L, L]$. The initial-value (Cauchy) problem is considered, and the initial conditions are written as
\begin{align}
	&u(x,0) = F_{1}(x), \quad u_{t} (x,0) = V_{1}(x), \label{uIC} \\
	&w(x,0) = F_{2}(x), \quad  w_{t} (x,0) = V_{2}(x). \label{wIC}
\end{align}
Firstly, we integrate \eqref{ueq} - \eqref{weq} in $x$ over the period $2L$ to obtain evolution equations of the form
\begin{align}
	\diffn{2}{ }{t} \int_{-L}^{L} u(x,t) \dd{x} + \varepsilon \delta \int_{-L}^{L} \lb u(x,t) - w(x,t) \rb \dd{x} &= 0, \label{ueqmean} \\
	\diffn{2}{ }{t} \int_{-L}^{L} w(x,t) \dd{x} - \varepsilon \gamma \int_{-L}^{L} \lb u(x,t) - w(x,t) \rb \dd{x} &= 0. \label{weqmean}
\end{align}
Denoting the mean value of $u$ and $w$ as
\begin{equation}
\langle u \rangle (t) := \frac{1}{2L} \int_{-L}^{L} u(x,t) \dd{x}, \quad \langle w \rangle (t) := \frac{1}{2L} \int_{-L}^{L} w(x,t) \dd{x},
\label{MeanVal}
\end{equation}
we can solve \eqref{ueqmean} - \eqref{weqmean} to obtain
\begin{align}
	\langle u \rangle &= d_{1} + \delta d_{2} \cos{\omega t} + d_{3} t + \delta d_{4} \sin{\omega t}, \label{uMeanVal} \\
	\langle w \rangle &= d_{1} - \gamma d_{2} \cos{\omega t} + d_{3} t - \gamma d_{4} \sin{\omega t}, \label{wMeanVal}
\end{align}
where $\omega = \sqrt{\varepsilon \lb \delta + \gamma \rb}$. Using the initial conditions \eqref{uIC} - \eqref{wIC} we can determine the values of the coefficients as
\begin{equation}
	d_{1} = \frac{\gamma \langle F_{1} \rangle + \delta \langle F_{2} \rangle}{\delta + \gamma}, \quad d_{2} = \frac{\langle F_{1} \rangle - \langle F_{2} \rangle}{\delta + \gamma}, \quad d_{3} = \frac{\gamma \langle V_{1} \rangle + \delta \langle V_{2} \rangle}{\omega \lb \delta + \gamma \rb}, \quad d_{4} = \frac{\langle V_{1} \rangle - \langle V_{2} \rangle}{\omega \lb \delta + \gamma \rb},
	\label{MeanValdi}
\end{equation}
and we have
\begin{equation}
	\langle F_{i}  \rangle= \int_{-L}^{L} F_{i}(x) \dd{x}, \quad \langle V_{i} \rangle = \int_{-L}^{L} V_{i}(x) \dd{x}, \quad i = 1,2.
	\label{MeanValIC}
\end{equation}
To simplify the problem we will consider initial conditions that satisfy the condition $d_{3} = d_{4} = 0$, that is
\begin{equation}
	\frac{1}{2L} \int_{-L}^{L} V_{i} \dd{x} = 0, \quad i = 1,2.
	\label{UWCondition2}
\end{equation}
This condition appears naturally in many physical applications and is imposed here to simplify our derivations, however, as was shown for the BKG equation, it can be relaxed (see the Appendix in \cite{Khusnutdinova19}). 

We subtract \eqref{uMeanVal} from $u$ and \eqref{wMeanVal} from $w$ to construct an equation with zero mean value, so we introduce $\tilde{u} = u - \langle u \rangle $ and $\tilde{w} = w - \langle w \rangle $ to obtain the problem for deviations
\begin{align}
	&\tilde{u}_{tt} - \tilde{u}_{xx} = \varepsilon \lsq \frac{1}{2} \lb \tilde{u}^2 \rb_{xx} + \langle u \rangle \tilde{u}_{xx} + \tilde{u}_{ttxx} - \delta \lb \tilde{u} - \tilde{w} \rb \rsq, \label{ueqtilde} \\
	&\tilde{w}_{tt} - c^2 \tilde{w}_{xx} = \varepsilon \lsq \frac{\alpha}{2} \lb \tilde{w}^2 \rb_{xx} + \alpha \langle w \rangle \tilde{w}_{xx} + \beta \tilde{w}_{ttxx} + \gamma \lb \tilde{u} - \tilde{w} \rb \rsq, \label{weqtilde}
\end{align}
where the expression for $\langle u \rangle$ and $\langle w \rangle$ can be found in \eqref{uMeanVal} and \eqref{wMeanVal} respectively. Note that this problem has variable coefficients, and the traditional procedure used to derive uni-directional models of the Korteweg-de Vries/Ostrovsky type is no longer applicable. The initial conditions become
\begin{align}
	&\tilde{u}(x,0) = \tilde{F}_{1}(x) = F_{1}(x) - \langle F_1 \rangle, \quad \tilde{u}_{t}(x,0) = 
	%\tilde{V}_{1}(x) = 
	V_{1}(x), \label{IC1} \\
	&\tilde{w}(x,0) = \tilde{F}_{2}(x) = F_{2}(x) - \langle F_2 \rangle, \quad \tilde{w}_{t}(x,0) = 
	%\tilde{V}_{2}(x) = 
	V_{2}(x), \label{IC2}
\end{align}
and, by construction, have zero mean value. The case of $c - 1 = \O{\varepsilon}$ was considered in \cite{Khusnutdinova19Book} and so we will only present the derivation for $c - 1 = \O{1}$. We will compare the cases in the results section using the derivation in \cite{Khusnutdinova19Book}.

In the case when $c - 1 = \O{1}$ the characteristic variables cannot be the same in each layer, and instead we have two distinct pairs of characteristic variables. In what follows we omit tildes and look for a weakly-nonlinear solution of the initial-value problem \eqref{ueqtilde} - \eqref{IC2} of the form
\begin{align}
	u(x,t) &= f_{1}^{-} \lb \xi_{-}, \tau, T \rb + f_{1}^{+} \lb \xi_{+}, \tau, T \rb + \sqrt{\varepsilon} P_{1} \lb \xi_{-}, \xi_{+}, \tau, T \rb + \varepsilon Q_{1} \lb \xi_{-}, \xi_{+}, \tau, T \rb + \varepsilon^{\frac{3}{2}} R_{1} \lb \xi_{-}, \xi_{+}, \tau, T \rb \notag \\
	&~~~+ \varepsilon^2 S_{1} \lb \xi_{-}, \xi_{+}, \tau, T \rb + \O{\varepsilon^{\frac{5}{2}}}, \label{uWNLc1} \\
	w(x,t) &= f_{2}^{-} \lb \nu_{-}, \tau, T \rb + f_{2}^{+} \lb \nu_{+}, \tau, T \rb + \sqrt{\varepsilon} P_{2} \lb \nu_{-}, \nu_{+}, \tau, T \rb + \varepsilon Q_{2} \lb \nu_{-}, \nu_{+}, \tau, T \rb + \varepsilon^{\frac{3}{2}} R_{2} \lb \nu_{-}, \nu_{+}, \tau, T \rb \notag \\
	&~~~+ \varepsilon^2 S_{2} \lb \nu_{-}, \nu_{+}, \tau, T \rb + \O{\varepsilon^{\frac{5}{2}}}, \label{wWNLc1}
\end{align}
where we use the characteristic and slow time variables
\begin{equation*}
	\xi_{\pm} = x \pm t, \quad \nu_{\pm} = x \pm ct, \quad \tau = \sqrt{\varepsilon} t, \quad T = \varepsilon t.
\end{equation*}
Note that the second time scale, $\tau$, arises from the terms $\langle u \rangle$ and $\langle w \rangle$ as we have the function $\cos(\omega t)$. We can extract the factor of $\sqrt{\varepsilon}$ from $\omega$ and write $\cos(\tilde{\omega}\tau)$, where $\tilde{\omega} = \sqrt{\delta + \gamma}$.

As we are considering the solution on the periodic domain, $u$ and $w$ are $2L$-periodic functions in $x$. Therefore we require that $f_{1}^{-}$ and $f_{1}^{+}$ are periodic in $\xi_{-}$ and $\xi_{+}$ respectively, and similarly $f_{2}^{-}$ and $f_{2}^{+}$ are periodic in $\nu_{-}$ and $\nu_{+}$ respectively. We then ensure that all functions in the constructed expansion are periodic as well. We will also construct our solution so that functions at each order have zero mean. In particular,
\begin{equation}
	\frac{1}{2L} \int_{-L}^{L} f_{1}^{\pm} \dd{\xi_{\pm}} = 0, \quad \frac{1}{2L} \int_{-L}^{L} f_{2}^{\pm} \dd{\nu_{\pm}} = 0.
	\label{fZM}
\end{equation}
To do this, for the functions in the expansion we will choose appropriate initial conditions, using those found earlier for the cRB equations, so that the initial conditions have zero mean. Where an integration occurs, we will choose any integration constants as the subtraction of the mean value of the function, to enforce zero mean conditions. This will be explained at each step where it occurs. This is consistent with the approach taken in \cite{Khusnutdinova19} where we worked in the space of functions with zero mean. In our previous work the integration constants were not explicitly stated in the derivation, but they were implemented to maintain zero mean conditions in the comparisons with the numerical solutions.

As mentioned earlier, the case when $c - 1 = \O{\varepsilon}$ was considered in \cite{Khusnutdinova19Book} and the expansions were in the same set of characteristic variables, while here we have two sets of characteristic variables. This means that, at some stage in our derivation, we may encounter a situation where we have a function of $\nu_{\pm}$ in the equation where the natural set of characteristic variables is $\xi_{\pm}$. In this case, we can rewrite one set of characteristic variables as a linear combination of the other set, so $\nu_{\pm}$ in terms of $\xi_{\pm}$, and then proceed. In each case, we will be looking at functions that are fully defined at an earlier order, but are being evaluated in terms of different characteristic variables. A similar situation was encountered in \cite{Khusnutdinova11} for the problem considered on the infinite domain, and we will introduce the linear combination when it is required.

We substitute \eqref{uWNLc1} and \eqref{wWNLc1} into \eqref{ueqtilde} and \eqref{weqtilde}, then compare at increasing powers of $\sqrt{\varepsilon}$. The equations are satisfied at leading order so we move onto terms at $\O{\sqrt{\varepsilon}}$. At this order we have
\begin{equation}
	-4 P_{1 \xi_{-} \xi_{+}} - 2 f^{-}_{1 \xi_{-} \tau} + 2 f^{+}_{1 \xi_{+} \tau} = 0.
	\label{P1eqpreavc1}
\end{equation}
We average \eqref{P1eqpreavc1} with respect to the fast spatial variable $x$ at constant $\xi_{-}$ and $\xi_{+}$ (see Ref.~\cite{Khusnutdinova19}). Averaging $P_{1 \xi_{-} \xi_{+}}$ at constant $\xi_{-}$ gives
\begin{equation}
	\frac{1}{2L} \int_{-L}^{L} P_{1 \xi_{-} \xi_{+}} \dd{x} = \frac{1}{4L} \int_{-2L - \xi_{-}}^{2L - \xi_{-}} P_{1 \xi_{-} \xi_{+}} \dd{\xi_{+}} = \frac{1}{4L} \lsq P_{1 \xi_{-}} \rsq_{-2L - \xi_{-}}^{2L - \xi_{-}} = 0,
	\label{Pav}
\end{equation}
with a similar result for averaging at constant $\xi_{+}$. Applying the averaging to \eqref{P1eqpreavc1} and requiring that $f_1^{\pm}$ have zero mean, we find
\begin{equation}
	f^{-}_{1 \xi_{-} \tau} = 0 \implies f^{-}_{1} = f^{-}_{1} \lb \xi_{-}, T \rb, \quad 
	f^{+}_{1 \xi_{+} \tau} = 0 \implies f^{+}_{1} = f^{+}_{1} \lb \xi_{+}, T \rb.
	\label{f1eqc1}
\end{equation}
A similar approach for the equation in $w$ gives
\begin{equation}
	f^{-}_{2} = f^{-}_{2} \lb \nu_{-}, T \rb \eqtext{and} f^{+}_{2} = f^{+}_{2} \lb \nu_{+}, T \rb.
	\label{f2eqc1}
\end{equation}
Substituting \eqref{f1eqc1} into \eqref{P1eqpreavc1} we obtain
\begin{equation}
	P_{1 \xi_{-} \xi_{+}} = 0 \implies P_{1} = g_{1}^{-} \lb \xi_{-}, \tau, T \rb + g_{1}^{+} \lb \xi_{+}, \tau, T \rb,
	\label{P1eqc1}
\end{equation}
and similarly for $w$ we obtain
\begin{equation}
	P_{2} = g_{2}^{-} \lb \nu_{-}, \tau, T \rb + g_{2}^{+} \lb \nu_{+}, \tau, T \rb.
	\label{P2eqc1}
\end{equation}
The initial conditions for $f_{1,2}^{\pm}$ are found by substituting \eqref{uWNLc1} into \eqref{uIC}, and also \eqref{wWNLc1} into \eqref{wIC}, then comparing terms at $\O{1}$, to obtain
\begin{align}
	f_{1}^{\pm}|_{T = 0} &= \frac{1}{2} \lb F_{1} \lb \xi_{\pm} \rb \pm \lb \int_{-L}^{\xi_{\pm}} V_{1} \lb \sigma \rb \dd{\sigma} - \tilde{V}_1^{\pm} \rb \rb, \label{f1ICc1} \\
	f_{2}^{\pm}|_{T = 0} &= \frac{1}{2c} \lb c F_{2} \lb \nu_{\pm} \rb \pm \lb \int_{-L}^{\nu_{\pm}} V_{2} \lb \sigma \rb \dd{\sigma} - \tilde{V}_2^{\pm} \rb \rb.
	\label{f2ICc1}
\end{align}
The values $\tilde{V}_1^{\pm}$ and $\tilde{V}_2^{\pm}$ are chosen in such a way to ensure that the initial conditions have zero mean. Explicitly, we calculate $\tilde{V}_1^{\pm}$ and $\tilde{V}_2^{\pm}$ as
\begin{equation*}
	\tilde{V}_1^{\pm} = \int_{-L}^{L} \int_{-L}^{\xi_{\pm}} V_1(\sigma) \dd{\sigma} \dd{\xi_{\pm}}, \quad \tilde{V}_2^{\pm} = \int_{-L}^{L} \int_{-L}^{\nu_{\pm}} V_2(\sigma) \dd{\sigma} \dd{\nu_{\pm}}.
\end{equation*}
This approach of enforcing zero mean will be used at various stages throughout our derivation. To find equations for $f_{1}^{\pm}$ and $f_{2}^{\pm}$, we next compare terms at $\O{\varepsilon}$, using the results from the previous order. Therefore, we find
\begin{align}
	-4Q_{1 \xi_{-} \xi_{+}} &= \lb 2 f_{1 T}^{-} + f_{1}^{-} f_{1 \xi_{-}}^{-} + d_{1} f_{1 \xi_{-}}^{-} + f_{1 \xi_{-} \xi_{-} \xi_{-}}^{-} \rb_{\xi_{-}} + \lb -2 f_{1 T}^{+} + f_{1}^{+} f_{1 \xi_{+}}^{+} + d_{1} f_{1 \xi_{+}}^{+} + f_{1 \xi_{+} \xi_{+} \xi_{+}}^{+} \rb_{\xi_{+}} \notag \\
	&~~~- \delta \lb f_{1}^{-} - f_{2}^{-} \rb - \delta \lb f_{1}^{+} - f_{2}^{+} \rb + 2 g_{1 \xi_{-} \tau}^{-} - 2 g_{1 \xi_{+} \tau}^{+} + d_{2} \delta \cos{\lb \tilde{\omega} \tau \rb} \lb f_{1 \xi_{-} \xi_{-}}^{-} + f_{1 \xi_{+} \xi_{+}}^{+} \rb \notag \\
	&~~~+ f_{1 \xi_{-} \xi_{-}}^{-} f_{1}^{+} + 2 f_{1 \xi_{-}}^{-} f_{1 \xi_{+}}^{+} + f_{1}^{-} f_{1 \xi_{+} \xi_{+}}^{+},
	\label{Q1eqpreavc1} \\
	-4 c^2 Q_{2 \nu_{-} \nu_{+}} &= \lb 2 c f_{2 T}^{-} + \alpha f_{2}^{-} f_{2 \nu_{-}}^{-} + \alpha d_{1} f_{2 \nu_{-}}^{-} + \beta c^2 f_{2 \nu_{-} \nu_{-} \nu_{-}}^{-} \rb_{\nu_{-}} + \gamma \lb f_{1}^{-} - f_{2}^{-} \rb \notag \\
	&~~~+ \lb -2 c f_{2 T}^{+} + \alpha f_{2}^{+} f_{2 \nu_{+}}^{+} + \alpha d_{1} f_{2 \nu_{+}}^{+} + \beta c^2 f_{2 \nu_{+} \nu_{+} \nu_{+}}^{+} \rb_{\nu_{+}} + \gamma \lb f_{1}^{+} - f_{2}^{+} \rb + 2c g_{2 \nu_{-} \tau}^{-} - 2 c g_{2 \nu_{+} \tau}^{+} \notag \\
	&~~~- \alpha d_{2} \gamma \cos{\lb \tilde{\omega} \tau \rb} \lb f_{2 \nu_{-} \nu_{-}}^{-} + f_{2 \nu_{+} \nu_{+}}^{+} \rb + \alpha \lb f_{2 \nu_{-} \nu_{-}}^{-} f_{2}^{+} + 2 f_{2 \nu_{-}}^{-} f_{2 \nu_{+}}^{+} + f_{2}^{-} f_{2 \nu_{+} \nu_{+}}^{+} \rb.
	\label{Q2eqpreavc1}
\end{align}
Averaging \eqref{Q1eqpreavc1} at constant $\xi_{-}$ or constant $\xi_{+}$ gives
\begin{equation}
	\pm 2 g_{1 \xi_{\pm} \tau} = d_{2} \delta \cos{\lb \tilde{\omega} \tau \rb} f_{1 \xi_{\pm} \xi_{\pm}}^{\pm} + A_{1} \lb \xi_{\pm}, T \rb,
	\label{g1diffc1}
\end{equation}
where
\begin{equation}
	A_{1} = \lb \mp 2 f_{1 T}^{\pm} + f_{1}^{\pm} f_{1 \xi_{\pm}}^{\pm} + d_{1} f_{1 \xi_{\pm}}^{\pm} + f_{1 \xi_{\pm} \xi_{\pm} \xi_{\pm}}^{\pm} \rb_{\xi_{\pm}} - \delta f_{1}^{\pm}.
	\label{A1c1}
\end{equation}
To avoid secular terms we require that $A_{1} = 0$. Therefore we obtain an Ostrovsky equation for $f_{1}^{\pm}$ of the form
\begin{equation}
	\lb \mp 2 f_{1 T}^{\pm} + f_{1}^{\pm} f_{1 \xi_{\pm}}^{\pm} + d_{1} f_{1 \xi_{\pm}}^{\pm} + f_{1 \xi_{\pm} \xi_{\pm} \xi_{\pm}}^{\pm} \rb_{\xi_{\pm}} = \delta f_{1}^{\pm}.
	\label{f1Ostc1}
\end{equation}
%{\color{red}
The Ostrovsky equation necessarily requires its regular solutions to have zero mean, and therefore as the initial condition also satisfies this condition, the function $f_{1}^{\pm}$ has zero mean. A similar approach for $w$, averaging at constant $\nu_{-}$ or constant $\nu_{+}$, leads to
\begin{equation}
	\pm 2 c g_{2 \nu_{\pm} \tau} = - d_{2} \gamma \cos{\lb \tilde{\omega} \tau \rb} f_{2 \nu_{\pm} \nu_{\pm}}^{\pm} + A_{2} \lb \nu_{\pm}, T \rb, 
	\label{g2diffc1}
\end{equation}
where again requiring $A_2 = 0$ to avoid secular terms gives an Ostrovsky equation for $f_{2}^{\pm}$ of the form
\begin{equation}
	\lb \mp 2 c f_{2 T}^{\pm} + \alpha f_{2}^{\pm} f_{2 \nu_{\pm}}^{\pm} + \alpha d_{1} f_{2 \nu_{\pm}}^{\pm} + \beta c^2 f_{2 \nu_{\pm} \nu_{\pm} \nu_{\pm}}^{\pm} \rb_{\nu_{\pm}} = \gamma f_{2}^{\pm},
	\label{f2Ostc1}
\end{equation}
and $f_{2}^{\pm}$ has zero mean by construction. Integrating \eqref{g1diffc1} we find an equation for $g_{1}^{\pm}$ of the form
\begin{equation}
	g_{1}^{\pm} = \pm \theta_{1} f_{1 \xi_{\pm}}^{\pm} + G_{1}^{\pm} \lb \xi_{\pm}, T \rb,
	\label{g1eqc1}
\end{equation}
and similarly for $g_2^{\pm}$ we can find
\begin{equation}
	g_{2}^{\pm} = \mp \theta_{2} f_{2 \nu_{\pm}}^{\pm} + G_{2}^{\pm} \lb \nu_{\pm}, T \rb,
	\label{g2eqc1}
\end{equation}
where $G_{1}^{\pm}$ and $G_{2}^{\pm}$ are functions to be found, and we introduce
\begin{equation}
	\theta_{1} = \frac{d_{2} \delta}{2 \tilde{\omega}} \sin{\lb \tilde{\omega} \tau \rb}, \quad \theta_{2} = \frac{\alpha d_{2} \gamma}{2c \tilde{\omega}} \sin{\lb \tilde{\omega} \tau \rb}.
	\label{theta12c1}
\end{equation}
As $f_{1,2}^{\pm}$ are zero-mean, if the functions $G_{1,2}^{\pm}$ are also zero-mean then so too is $g_{1,2}^{\pm}$. We will find an equation for $G_{1,2}^{\pm}$ at the next order.
Substituting \eqref{f1Ostc1} and \eqref{g1eqc1} into \eqref{Q1eqpreavc1} and integrating we obtain
\begin{equation}
	Q_{1} = h_{1}^{-} \lb \xi_{-}, \tau, T \rb + h_{1}^{+} \lb \xi_{+}, \tau, T \rb + h_{1c} \lb \xi_{-}, \xi_{+}, T \rb + \hat{f}_{2}^{-} \lb \nu_{-}, T \rb + \hat{f}_{2}^{+} \lb \nu_{+}, T \rb,
	\label{Q1eqc1}
\end{equation}
where
\begin{equation}
	h_{1c} = - \frac{1}{4} \lb f_{1 \xi_{-}}^{-} \smashoperator{\int_{-L}^{\xi_{+}}}\!\! f_{1}^{+} \!\! \lb \sigma \rb\! \dd{\sigma} + 2 f_{1}^{-} f_{1}^{+} + f_{1 \xi_{+}}^{+} \smashoperator{\int_{-L}^{\xi_{-}}}\!\! f_{1}^{-}\!\! \lb \sigma \rb\! \dd{\sigma} \rb - \tilde{h}_{1c},
	\label{h1cc1}
\end{equation}
and $\tilde{h}_{1c}$ is the mean value of the function, so that the function $h_{1c}$ has zero mean. The function $\hat{f}_{2}^{\pm} \lb \nu_{\pm}, T \rb$ is a function of $\nu_{\pm}$, not $\xi_{\pm}$, however we can write the characteristic variables $\nu_{\pm}$ in terms of $\xi_{\pm}$ as
\begin{equation}
	\nu_{-} = \frac{(1 + c)\xi_{-} + (1 - c)\xi_{+}}{2}, \quad \nu_{+} = \frac{(1 - c)\xi_{-} + (1 + c)\xi_{+}}{2}.
	\label{LinearChar}
\end{equation}
Therefore, we can rewrite our functions $f_{1,2}^{\pm}$ in terms of the other set of characteristic variables, and so we find $\hat{f}_{2}^{\pm} \lb \nu_{\pm}, T \rb$ can be written as
\begin{equation}
	 \hat{f}_{2}^{\pm} = \frac{\delta}{c^2 - 1} \int_{-L}^{\nu_{\pm}} \int_{-L}^{z} f_{2}^{\pm} \lb y, T \rb \dd{y} \dd{z} - c_{\hat{f}_2},
	\label{f2hatc1}
\end{equation}
where $c_{\hat{f}_{2}}$ is the appropriately chosen integration constant.  Similarly, substituting \eqref{f2Ostc1} and \eqref{g2eqc1} into \eqref{Q2eqpreavc1} and integrating we find
\begin{equation}
	Q_{2} = h_{2}^{-} \lb \nu_{-}, \tau, T \rb + h_{2}^{+} \lb \nu_{+}, \tau, T \rb + h_{2c} \lb \nu_{-}, \nu_{+}, T \rb + \hat{f}_{1}^{-} \lb \xi_{-}, T \rb + \hat{f}_{1}^{+} \lb \xi_{+}, T \rb,
	\label{Q2eqc1}
\end{equation}
where
\begin{equation}
	h_{2c} = - \frac{\alpha}{4c^2} \lb f_{2 \nu_{-}}^{-} \smashoperator{\int_{-L}^{\nu_{+}}}\!\! f_{2}^{+} \!\! \lb \sigma \rb\! \dd{\sigma} + 2 f_{2}^{-} f_{2}^{+} + f_{2 \nu_{+}}^{+} \smashoperator{\int_{-L}^{\nu_{-}}}\!\! f_{2}^{-}\!\! \lb \sigma \rb\! \dd{\sigma} \rb - \tilde{h}_{2c},
	\label{h2cc1}
\end{equation}
and $h_{2c}$ is the appropriately chosen integration constant. Rewriting $\xi_{\pm}$ in terms of $\nu_{\pm}$ in a similar way to \eqref{LinearChar}, we have
\begin{equation}
	\hat{f}_{1}^{\pm} = - \frac{\gamma}{c^2 - 1} \int_{-L}^{\xi_{\pm}} \int_{-L}^{z} f_{1}^{\pm} \lb y, T \rb \dd{y} \dd{z} - c_{\hat{f}_{1}},
	\label{f1hatc1}
\end{equation}
with $c_{\hat{f}_{1}}$ an appropriately chosen integration constant so that $\hat{f}_{1}^{\pm}$ has zero mean.
As before we find the initial condition for $G_{1,2}^{\pm}$ by substituting \eqref{uWNLc1} and \eqref{wWNLc1} into \eqref{uIC} and \eqref{wIC}, comparing terms at $\O{\sqrt{\varepsilon}}$. Therefore, taking account of the results found for \eqref{g1eqc1} and \eqref{g2eqc1}, and noting that $\theta_{1}|_{T = 0} = 0$, we obtain
\begin{equation}
	\left. G_{1}^{\pm} \right|_{T = 0} = 0 \quad \text{and} \quad \left. G_{2}^{\pm} \right|_{T = 0} = 0.
	\label{GICc1}
\end{equation}
We now aim to find equations for $G_{1,2}^{\pm}$ and $h_{1,2}^{\pm}$ at the next order. Comparing terms at $\O{\varepsilon^{3/2}}$ in the expansion of \eqref{ueqtilde} and \eqref{weqtilde}, taking account of results at previous orders, we have
\begin{align}
	-4R_{1 \xi_{-} \xi_{+}} &= \lb 2 g_{1 T}^{-} + \lb f_{1}^{-} g_{1}^{-} \rb_{\xi_{-}} + d_{1} g_{1 \xi_{-}}^{-} + g_{1 \xi_{-} \xi_{-} \xi_{-}}^{-} \rb_{\xi_{-}} + \lb - 2 g_{1 T}^{+} + \lb f_{1}^{+} g_{1}^{+} \rb_{\xi_{+}} + d_{1} g_{1 \xi_{+}}^{+} + g_{1 \xi_{+} \xi_{+} \xi_{+}}^{+} \rb_{\xi_{+}} \notag \\
	&~~~- \delta \lb g_{1}^{-} - g_{2}^{-} \rb - \delta \lb g_{1}^{+} - g_{2}^{+} \rb + 2 h_{1 \xi_{-} \tau}^{-} - 2 h_{1 \xi_{+} \tau}^{+} - g_{1 \tau \tau}^{-} - g_{1 \tau \tau}^{+} + d_{2} \delta \cos{\lb \tilde{\omega} \tau \rb} \lb g_{1 \xi_{-} \xi_{-}}^{-} + g_{1 \xi_{+} \xi_{+}}^{+} \rb \notag \\
	&~~~+ g_{1 \xi_{-} \xi_{-}}^{-} f_{1}^{+} + 2 g_{1 \xi_{-}}^{-} f_{1 \xi_{+}}^{+} + g_{1}^{-} f_{1 \xi_{+} \xi_{+}}^{+} + g_{1 \xi_{+} \xi_{+}}^{+} f_{1}^{-} + 2 g_{1 \xi_{+}}^{+} f_{1 \xi_{-}}^{-} + g_{1}^{+} f_{1 \xi_{-} \xi_{-}}^{-}, \label{R1eqpreavc1} \\
	-4 c^2 R_{2 \nu_{-} \nu_{+}} &= \lb 2 c g_{2 T}^{-} + \alpha \lb f_{2}^{-} g_{2}^{-} \rb_{\nu_{-}} + \alpha d_{1} g_{2 \nu_{-}}^{-} + \beta c^2 g_{2 \nu_{-} \nu_{-} \nu_{-}}^{-} \rb_{\nu_{-}} + \gamma \lb g_{1}^{-} - g_{2}^{-} \rb \notag \\
	&~~~+ \lb - 2 c g_{2 T}^{+} + \alpha \lb f_{2}^{+} g_{2}^{+} \rb_{\nu_{+}} + \alpha d_{1} g_{2 \nu_{+}}^{+} + \beta c^2 g_{2 \nu_{+} \nu_{+} \nu_{+}}^{+} \rb_{\nu_{+}} + \gamma \lb g_{1}^{+} - g_{2}^{+} \rb \notag \\
	&~~~+ 2 c h_{2 \nu_{-} \tau}^{-} - 2 c h_{2 \nu_{+} \tau}^{+} - g_{2 \tau \tau}^{-} - g_{2 \tau \tau}^{+} - \alpha d_{2} \gamma \cos{\lb \tilde{\omega} \tau \rb} \lb g_{2 \nu_{-} \nu_{-}}^{-} + g_{2 \nu_{+} \nu_{+}}^{+} \rb \notag \\
	&~~~+ \alpha \lsq g_{2 \nu_{-} \nu_{-}}^{-} f_{2}^{+} + 2 g_{2 \nu_{-}}^{-} f_{2 \nu_{+}}^{+} + g_{2}^{-} f_{2 \nu_{+} \nu_{+}}^{+} + g_{2 \nu_{+} \nu_{+}}^{+} f_{2}^{-} + 2 g_{2 \nu_{+}}^{+} f_{2 \nu_{-}}^{-} + g_{2}^{+} f_{2 \nu_{-} \nu_{-}}^{-} \rsq. \label{R2eqpreavc1}
\end{align}
Substituting \eqref{g1eqc1} into \eqref{R1eqpreavc1} and averaging at constant $\xi_{-}$ or constant $\xi_{+}$ leads to
\begin{align}
	\pm 2 h_{1 \xi_{\pm} \tau}^{\pm} &= \pm \theta_{1} \lb \mp 2 f_{1 T}^{\pm} + f_{1}^{\pm} f_{1 \xi_{\pm}}^{\pm} + d_{1} f_{1 \xi_{\pm}}^{\pm} + f_{1 \xi_{\pm} \xi_{\pm} \xi_{\pm}} \rb_{\xi_{\pm} \xi_{\pm}} + \lb \mp 2 G_{1 T}^{\pm} + \lb f_{1}^{\pm} G_{1}^{\pm} \rb_{\xi_{\pm}} + d_{1} G_{1 \xi_{\pm}}^{\pm} + G_{1 \xi_{\pm} \xi_{\pm} \xi_{\pm}} \rb_{\xi_{\pm}} \notag \\
	&~~~\mp \theta_{1} \delta f_{1}^{\pm} - \delta G_{1 \xi_{\pm}}^{\pm} \pm \theta_{1} \tilde{\omega}^2 f_{1 \xi_{\pm}}^{\pm} \pm \theta_{1} d_{2} \delta \cos{\lb \tilde{\omega} \tau \rb} f_{1 \xi_{\pm} \xi_{\pm} \xi_{\pm}}^{\pm}. 
	\label{R1averagedc1}
\end{align}
If we differentiate \eqref{f1Ostc1} with respect to the appropriate characteristic variable, we can eliminate some terms from \eqref{R1averagedc1} to obtain an expression for $h_{1 \xi_{\pm} \tau}^{\pm}$ of the form
\begin{equation}
	2 h_{1 \xi_{\pm} \tau}^{\pm} = \theta_1 \tilde{\omega}^2 f_{1 \xi_{\pm}}^{\pm} + \theta_{1} d_{2} \delta \cos{\lb \tilde{\omega} \tau \rb} f_{1 \xi_{\pm} \xi_{\pm} \xi_{\pm}}^{\pm} + \tilde{G}_{1}^{\pm} \lb \xi_{\pm}, T \rb,
	\label{h1diffc1}
\end{equation}
where
\begin{equation}
	\tilde{G}_{1}^{\pm} = \lb \mp 2 G_{1 T}^{\pm} + \lb f_{1}^{\pm} G_{1}^{\pm} \rb_{\xi_{\pm}} + d_{1} G_{1 \xi_{\pm}}^{\pm} + G_{1 \xi_{\pm} \xi_{\pm} \xi_{\pm}}^{\pm} \rb_{\xi_{\pm}} - \delta G_{1}^{\pm}.
	\label{G1tildeeqc1}
\end{equation}
To avoid secular terms we require that $\tilde{G}_{1}^{\pm} = 0$ and therefore we have an equation for $G_{1}^{\pm}$ of the form
\begin{equation}
	\lb \mp 2 G_{1 T}^{\pm} + \lb f_{1}^{\pm} G_{1}^{\pm} \rb_{\xi_{\pm}} + d_{1} G_{1 \xi_{\pm}}^{\pm} + G_{1 \xi_{\pm} \xi_{\pm} \xi_{\pm}}^{\pm} \rb_{\xi_{\pm}} = \delta G_{1}^{\pm}.
	\label{G1eqc1}
\end{equation}
Taking account of the initial condition in \eqref{GICc1} and the form of \eqref{G1eqc1}, we can clearly see that $G_{1}^{\pm} \equiv 0$ in this derivation and therefore in all subsequent steps we will omit all terms in $G_{1}^{\pm}$. Referring back to \eqref{g1eqc1}, we now see that $g_1^{\pm}$ has zero mean. Integrating \eqref{h1diffc1} we obtain
\begin{equation}
	h_{1}^{\pm} = -\frac{\tilde{\omega}^2 \rho_{1}}{2} f_{1}^{\pm} - \frac{\tilde{\omega}^2 \rho_{1}^2}{2} f_{1 \xi_{\pm} \xi_{\pm}}^{\pm} + \phi_{1}^{\pm} \lb \xi_{\pm}, T \rb,
	\label{h1eqc1}
\end{equation}
where the functions $\phi_{1}^{\pm}$ are to be found at the next order and
\begin{equation}
	\rho_{1} = \pinv_\tau {\theta_{1}} = \frac{d_{2} \delta}{2 \tilde{\omega}^2} \cos{\lb \tilde{\omega} \tau \rb}.
	\label{rho1c1}
\end{equation}
%{\color{red}
If the function $\phi_{1}^{\pm}$ is constructed to have zero mean, then so does $h_{1}^{\pm}$. Similarly, from the equation for $w$, we find 
\begin{equation}
	\pm 2 c h_{2 \xi_{\pm} \tau}^{\pm} = \mp \theta_2 \tilde{\omega}^2 f_{2 \nu_{\pm}}^{\pm} \pm \theta_{2} d_{2} \gamma \cos{\lb \tilde{\omega} \tau \rb} f_{2 \nu_{\pm} \nu_{\pm} \nu_{\pm}}^{\pm} + \tilde{G}_{2}^{\pm} \lb \nu_{\pm}, T \rb,
	\label{h2diffc1}
\end{equation}
where
\begin{equation}
	\tilde{G}_{2}^{\pm} \lb \nu_{\pm}, T \rb = \mp 2 c G_{2 T \nu_{\pm}}^{\pm} + \alpha \lb f_{2}^{\pm} G_{2}^{\pm} \rb_{\nu_{\pm} \nu_{\pm}} + \alpha d_{1} G_{2 \nu_{\pm} \nu_{\pm}}^{\pm} + \beta c^2 G_{2 \nu_{\pm} \nu_{\pm} \nu_{\pm} \nu_{\pm}}^{\pm} - \gamma G_{2}^{\pm}.
	\label{G2tildeeqc1}
\end{equation}
We require that $\tilde{G}_{2}^{\pm} = 0$ and, by the same argument as for $G_{1}^{\pm}$ we have $G_{2}^{\pm} \equiv 0$ and therefore omit it from all subsequent derivations, meaning $g_{2}^{\pm}$ has zero mean. Integrating \eqref{h2diffc1} we obtain
\begin{equation}
	h_{2}^{\pm} = \frac{\tilde{\omega}^2 \rho_ {2}}{2c} f_{2}^{\pm} - \frac{\tilde{\omega}^2 \rho_{2}^2}{2} f_{2 \nu_{\pm} \nu_{\pm}}^{\pm} + \phi_{2}^{\pm} \lb \nu_{\pm}, T \rb,
	\label{h2eqc1}
\end{equation}
where again we need to find the function $\phi_{2}^{\pm}$ and
\begin{equation}
	\rho_{2} = \pinv_\tau {\theta_{2}} = \frac{\alpha d_{2} \gamma}{2 c \tilde{\omega}^2} \cos{\lb \tilde{\omega} \tau \rb}.
	\label{rho2c1}
\end{equation}
If the function $\phi_{2}^{\pm}$ is constructed to have zero mean, then so does $h_{2}^{\pm}$.
%} 
Substituting \eqref{h1eqc1} into \eqref{R1eqpreavc1} and integrating with respect to the appropriate characteristic variables we find
\begin{equation}
	R_{1} = \psi_{1}^{-} \lb \xi_{-}, \tau, T \rb + \psi_{1}^{+} \lb \xi_{+}, \tau, T \rb + \psi_{1c} \lb \xi_{-}, \xi_{+}, T \rb + \hat{g}_{2}^{-} + \hat{g}_{2}^{+},
	\label{R1eqc1}
\end{equation}
where
\begin{equation}
	\psi_{1c} = -\frac{\theta_{1}}{4} \lsq f_{1 \xi_{+} \xi_{+}}^{+} \int_{-L}^{\xi_{-}} f_{1}^{-} \lb \sigma \rb \dd{\sigma} - f_{1}^{+} f_{1 \xi_{-}}^{-} + f_{1}^{-} f_{1 \xi_{+}}^{+} - f_{1 \xi_{-} \xi_{-}}^{-} \int_{-L}^{\xi_{+}} f_{1}^{+} \lb \sigma \rb \dd{\sigma} \rsq - c_{\psi_{1c}},
	\label{psi1cc1}
\end{equation}
with appropriate integration constant to maintain zero mean. The terms $\hat{g}_{2}^{\pm}$ can be found by replacing $f$ with $g$ in \eqref{f2hatc1}. 
Similarly, substituting \eqref{h2eqc1} into \eqref{R2eqpreavc1} and integrating gives
\begin{equation}
	R_{2} = \psi_{2}^{-} \lb \nu_{-}, \tau, T \rb + \psi_{2}^{+} \lb \nu_{+}, \tau, T \rb + \psi_{2c} \lb \xi_{-}, \xi_{+}, T \rb + \hat{g}_{1}^{-} + \hat{g}_{1}^{+},
	\label{R2eqc1}
\end{equation}
where
\begin{equation}
	\psi_{2c} = -\frac{\alpha \theta_{2}}{4 c^2} \lsq f_{2 \nu_{-} \nu_{-}}^{-} \int_{-L}^{\nu_{+}} f_{2}^{+} \lb \sigma \rb \dd{\sigma} - f_{2 \nu_{-}}^{-} f_{2}^{+} + f_{2}^{-} f_{2 \nu_{+}}^{+} - f_{2 \nu_{+} \nu_{+}}^{+} \int_{-L}^{\nu_{-}} f_{2}^{-} \lb \sigma \rb \dd{\sigma} \rsq - c_{\psi_{2c}},
	\label{psi2cc1}
\end{equation}
with appropriate integration constant to maintain zero mean. The expression for $\hat{g}_{1}^{\pm}$ can be found by replacing $f$ with $g$ in \eqref{f1hatc1}. The exact form of these terms is omitted as we are only interested in terms up to and including $\O{\varepsilon}$.

%, however in theory the derivation could be continued to any desired order. We have not included an integration constant in the expressions for $R_{1,2}$ as they are included within the constituent functions.

To find the initial condition for the function $\phi_{1}^{\pm}$, we substitute \eqref{uWNLc1} into \eqref{uIC} and comparing terms at $\O{\varepsilon}$, taking account of \eqref{h1eqc1}. Similarly, we can substitute \eqref{wWNLc1} into \eqref{wIC} and take account of \eqref{h2eqc1} to find an initial condition for $\phi_{2}^{\pm}$. Therefore 
\begin{align}
	\phi_{1}^{\pm} &= \frac{1}{2} \lb J_{1}^{\pm} \mp \int_{-L}^{\xi_{\pm}} K_{1} \lb \sigma \rb \dd{\sigma} \rb - \tilde{\phi}_1^{\pm}, \label{phi1ICc1} \\
	 \phi_{2}^{\pm} &= \frac{1}{2c} \lb J_{2}^{\pm} \mp \int_{-L}^{ \nu_{\pm}} K_{2} \lb \sigma \rb \dd{\sigma} \rb - \tilde{\phi}_{2}^{\pm},
	\label{phi2ICc1}
\end{align}
where
\begin{align}
	J_{1} &= \frac{\tilde{\omega}^2 \rho_{1}}{2} \lb f_{1}^{-} + f_{1}^{+} \rb + \frac{\tilde{\omega}^2 \rho_{1}^2}{2} \lb f_{1 \xi_{-} \xi_{-}}^{-} + f_{1 \xi_{+} \xi_{+}}^{+} \rb - h_{1c} - \frac{\delta}{c^2 - 1} \int_{-L}^{\nu_{-}} \int_{-L}^{v} f_{2}^{-} \lb u, T \rb \dd{u} \dd{v} \notag \\
	&~~~- \frac{\delta}{c^2 - 1} \int_{-L}^{\nu_{+}} \int_{-L}^{v} f_{2}^{+} \lb u, T \rb \dd{u} \dd{v}, \notag \\
	K_{1} &= f_{1 T}^{-} + f_{1 T}^{+} - \frac{\tilde{\omega}^2 \rho_{1}}{2} \lb f_{1 \xi_{-}}^{-} - f_{1 \xi_{+}}^{+} \rb + \frac{\tilde{\omega}^2 \rho_{1}^2}{2} \lb f_{1 \xi_{-} \xi_{-} \xi_{-}}^{-} - f_{1 \xi_{+} \xi_{+} \xi_{+}}^{+} \rb - h_{1c \xi_{-}} + h_{1c \xi_{+}} \notag \\
	&~~~- \frac{c \delta}{c^2 - 1} \int_{-L}^{\nu_{-}} f_{2}^{-} \lb y, T \rb \dd{y} + \frac{c \delta}{c^2 - 1} \int_{-L}^{\nu_{+}} f_{2}^{+} \lb y, T \rb \dd{y}, \label{JK1c1} \\
	J_{2} &= -\frac{\tilde{\omega}^2 \rho_{2}}{2c} \lb f_{2}^{-} + f_{2}^{+} \rb + \frac{\tilde{\omega}^2 \rho_{2}^2}{2} \lb f_{2 \nu_{-} \nu_{-}}^{-} + f_{2 \nu_{+} \nu_{+}}^{+} \rb  - h_{2c} + \frac{\gamma}{c^2 - 1} \int_{-L}^{\xi_{-}} \int_{-L}^{z} f_{1}^{-} \lb y, T \rb \dd{y} \dd{z} \notag \\
	&~~~+ \frac{\gamma}{c^2 - 1} \int_{-L}^{\xi_{+}} \int_{-L}^{z} f_{1}^{+} \lb y, T \rb \dd{y} \dd{z}, \notag \\
	K_{2} &= f_{2 T}^{-} + f_{2 T}^{+} + \frac{\tilde{\omega}^2 \rho_{2}}{2} \lb f_{2 \nu_{-}}^{-} - f_{2 \nu_{+}}^{+} \rb + \frac{c \tilde{\omega}^2 \rho_{2}^2}{2} \lb f_{2 \nu_{-} \nu_{-} \nu_{-}}^{-} - f_{2 \nu_{+} \nu_{+} \nu_{+}}^{+} \rb - c h_{2c \nu_{-}} + c h_{2c \nu_{+}} \notag \\
	&~~~+ \frac{\gamma}{c^2 - 1} \int_{-L}^{\xi_{-}} f_{1}^{-} \lb y, T \rb \dd{y} - \frac{\gamma}{c^2 - 1} \int_{-L}^{\xi_{+}} f_{1}^{+} \lb y, T \rb \dd{y}, \\
	\tilde{\phi}_1^{\pm} &= \int_{-L}^{L} \frac{1}{2} \lb J_{1}^{\pm} \mp \int_{-L}^{\xi_{\pm}} K_{1} \lb \sigma \rb \rb \dd{\xi_{\pm}}, \quad 
	\tilde{\phi}_{2}^{\pm} = \int_{-L}^{L} \frac{1}{2c} \lb J_{2}^{\pm} \mp \int_{-L}^{\nu_{\pm}} K_{2} \lb \sigma \rb \dd{\sigma} \rb \dd{\nu_{\pm}}.
\end{align}
We now have expressions up to $\O{\varepsilon}$, however we still need to find an equation governing $\phi_{1,2}^{\pm}$, therefore we compare terms at $\O{\varepsilon^2}$, taking account of all previous results. All coupling terms between left- and right-propagating waves are gathered in one function, and terms of the type of \eqref{f2hatc1} and \eqref{f1hatc1} are gathered in another function for convenience, as we do not require them to determine $\phi_{1,2}^{\pm}$. Gathering terms at $\O{\varepsilon^2}$ we have
\begin{align}
	- 4 S_{1 \xi_{-} \xi_{+}} &= - f_{1 T T}^{-} - f_{1 T T}^{+} - 2 g_{1 \tau T}^{-} - 2 g_{1 \tau T}^{+} - h_{1 \tau \tau}^{-} - h_{1 \tau \tau}^{+} + 2 h_{1 \xi_{-} T}^{-} - 2 h_{1 \xi_{+} T}^{+} + 2 \psi_{1 \xi_{-} \tau}^{-} - 2 \psi_{1 \xi_{+} \tau}^{+} \notag \\
	&~~~+ \lb f_{1}^{-} h_{1}^{-} \rb_{\xi_{-} \xi_{-}} + \lb f_{1}^{+} h_{1}^{+} \rb_{\xi_{+} \xi_{+}} + \frac{1}{2} \lb g_{1}^{-^{2}} \rb_{\xi_{-} \xi_{-}} + \frac{1}{2} \lb g_{1}^{+^{2}} \rb_{\xi_{+} \xi_{+}} + h_{1 \xi_{-} \xi_{-} \xi_{-} \xi_{-}}^{-} + h_{1 \xi_{+} \xi_{+} \xi_{+} \xi_{+}}^{+} \notag \\
	&~~~+ \lb d_{1} + d_{2} \delta \cos{\lb \tilde{\omega} \tau \rb} \rb \lb h_{1 \xi_{-} \xi_{-}}^{-} + h_{1 \xi_{+} \xi_{+}}^{+} \rb - 2 g_{1 \xi_{-} \xi_{-} \xi_{-} \tau}^{-} + 2 g_{1 \xi_{+} \xi_{+} \xi_{+} \tau}^{+} - 2 f_{1 \xi_{-} \xi_{-} \xi_{-} T}^{-} + 2 f_{1 \xi_{+} \xi_{+} \xi_{+} T}^{+} \notag \\
	&~~~- \delta \lb h_{1}^{-} - h_{2}^{-} + h_{1}^{+} - h_{2}^{+} \rb - 4 \mu_{1c} - 4 \Upsilon_{1},
	\label{S1eqpreavc1}
\end{align}
and
\begin{align}
	-4c^2 S_{2 \nu_{-} \nu_{+}} &= -f_{2 T T}^{-} - f_{2 T T}^{+} - 2 g_{2 \tau T}^{-} - 2 g_{2 \tau T}^{+} - h_{2 \tau \tau}^{-} - h_{2 \tau \tau}^{+} + 2 c h_{2 \nu_{-} T}^{-} - 2 c h_{2 \nu_{+} T}^{+} + 2 c \psi_{2 \nu_{-} \tau}^{-} - 2 c \psi_{2 \nu_{+} \tau}^{+} \notag \\
	&~~~+ \alpha \lb f_{2}^{-} h_{2}^{-} \rb_{\nu_{-} \nu_{-}} + \alpha \lb f_{2}^{+} h_{2}^{+} \rb_{\nu_{+} \nu_{+}} + \alpha \lb d_{1} - d_{2} \gamma \cos{\lb \tilde{\omega} \tau \rb} \rb \lb h_{2 \nu_{-} \nu_{-}}^{-} + h_{2 \nu_{+} \nu_{+}}^{+} \rb + \frac{\alpha}{2} \lb g_{2}^{-^{2}} \rb_{\nu_{-} \nu_{-}} \notag \\
	&~~~+ \frac{\alpha}{2} \lb g_{2}^{+^{2}} \rb_{\nu_{+} \nu_{+}} + \beta c^2 h_{2 \nu_{-} \nu_{-} \nu_{-} \nu_{-}}^{-} + \beta c^2 h_{2 \nu_{+} \nu_{+} \nu_{+} \nu_{+}}^{+} - 2 \beta c g_{2 \nu_{-} \nu_{-} \nu_{-} \tau}^{-} + 2 \beta c g_{2 \nu_{+} \nu_{+} \nu_{+} \tau}^{+} \notag \\
	&~~~- 2 \beta c f_{2 \nu_{-} \nu_{-} \nu_{-} T}^{-} + 2 \beta c f_{2 \nu_{+} \nu_{+} \nu_{+} T}^{+} + \gamma \lb h_{1}^{-} - h_{2}^{-} + h_{1}^{+} - h_{2}^{+} \rb - 4 \mu_{2c} - 4 \Upsilon_{2},
	\label{S2eqpreavc1}
\end{align}
where $\mu_{1c}$ and $\mu_{2c}$ are the coupling terms at this order, while $\Upsilon_{1}$ and $\Upsilon_{2}$ are the terms involving $\hat{f}_{2}$ or equivalent. We average \eqref{S1eqpreavc1} at constant $\xi_{-}$ or constant $\xi_{+}$, or average \eqref{S2eqpreavc1} at constant $\nu_{-}$ or constant $\nu_{+}$, integrate with respect to $\tau$ and rearrange to obtain
\begin{equation}
	2 \psi_{1 \xi_{\pm}} = H_{1}^{\pm} \lb \xi_{\pm}, \tau, T \rb + \hat{H}_{1}^{\pm} \lb \xi_{\pm}, T \rb \tau, \quad 2 c \psi_{2 \nu_{\pm}} = H_{2}^{\pm} \lb \nu_{\pm}, \tau, T \rb + \hat{H}_{2}^{\pm} \lb \nu_{\pm}, T \rb \tau,
	\label{psi12eqdiffc1}
\end{equation}
where the functions $H_{1,2}^{\pm}, \hat{H}_{1,2}^{\pm}$ can be found from \eqref{S1eqpreavc1} or \eqref{S2eqpreavc1}. To avoid secular terms we require that $\hat{H}_{1,2}^{\pm} = 0$ and this allows us to find equations for $\phi_{1,2}^{\pm}$. Therefore we look for terms in \eqref{S1eqpreavc1} that depend only on $\xi_{\pm}$ and $T$. Following this approach we obtain
\begin{align}
	\lb \mp 2 \phi_{1 T}^{\pm} + \lb f_{1}^{\pm} \phi_{1}^{\pm} \rb_{\xi_{\pm}} + d_{1} \phi_{1 \xi_{\pm}}^{\pm} + \phi_{1 \xi_{\pm} \xi_{\pm} \xi_{\pm}}^{\pm} \rb_{\xi_{\pm}} &= \delta \phi_{1}^{\pm} + f_{1 T T}^{\pm} \mp 2 f_{1 \xi_{\pm} \xi_{\pm} \xi_{\pm} T}^{\pm} \notag \\
	&~~~+ \frac{\tilde{\omega}^2 \tilde{\theta}_{1}^{2}}{2} f_{1 \xi_{\pm} \xi_{\pm}}^{\pm} - \frac{\tilde{\theta}_{1}^{2}}{2} \lb f_{1 \xi_{\pm}}^{\pm^{2}} \rb_{\xi_{\pm} \xi_{\pm}},
	\label{phi1eqc1}
\end{align}
and from \eqref{S2eqpreavc1} we have
\begin{align}
	\lb \mp 2c \phi_{2 T}^{\pm} + \alpha \lb f_{2}^{\pm} \phi_{2}^{\pm} \rb_{\nu_{\pm}} + \alpha d_{1} \phi_{2 \nu_{\pm}}^{\pm} + \beta c^2 \phi_{2 \nu_{\pm} \nu_{\pm} \nu_{\pm}}^{\pm} \rb_{\nu_{\pm}} &= \gamma \phi_{2}^{\pm} + f_{2 T T}^{\pm} \mp 2 c \beta f_{2 \nu_{\pm} \nu_{\pm} \nu_{\pm} T}^{\pm} \notag \\
	&~~~+ \frac{\tilde{\omega}^2 \tilde{\theta}_{2}^{2}}{2} f_{2 \nu_{\pm} \nu_{\pm}}^{\pm} - \frac{\alpha \tilde{\theta}_{2}^{2}}{2} \lb f_{2 \nu_{\pm}}^{\pm^{2}} \rb_{\nu_{\pm} \nu_{\pm}},
	\label{phi2eqc1}
\end{align}
where
\begin{equation}
	\tilde{\theta}_{1} = \frac{\theta_{1}}{\sin{\lb \tilde{\omega} \tau \rb}} = \frac{d_{2} \delta}{2 \tilde{\omega}}, \quad \tilde{\theta}_{2} = \frac{\theta_{2}}{\sin{\lb \tilde{\omega} \tau \rb}} = \frac{\alpha d_{2} \gamma}{2 c \tilde{\omega}}.
	\label{theta12tildec1}
\end{equation}
%}%
We have now defined all functions up to and including $\O{\varepsilon}$ and so stop our derivation, however in theory this could be continued to any order. We also note that every function in this expansion has been constructed to have zero mean, either directly or by the appropriate choice of integration constant.

% Results
\section{Validity of Weakly-Nonlinear Solution}
\label{sec:Validity}
In Section \ref{sec:WNL} we constructed the weakly-nonlinear solution to the original cRB equations \eqref{ueq} - \eqref{wIC}, when the characteristic speeds are essentially distinct (high-contrast case). We now confirm the validity of the constructed expansions by numerically solving the system \eqref{ueq} - \eqref{weq} and comparing this direct numerical solution to the constructed solution \eqref{uWNLc1} and \eqref{wWNLc1} with an increasing number of terms included. This was constructed in Section \ref{sec:WNL}, so we need to solve \eqref{f1eqc1}, \eqref{f2eqc1} for the leading-order solution and \eqref{phi1eqc1}, \eqref{phi2eqc1} for the solution up to and including terms at $\O{\varepsilon}$.

In order to numerically solve these equations, in this section and the subsequent sections, we implement three pseudo-spectral numerical schemes described in the appendices. For the coupled Boussinesq equations we use the methods in Appendix A and for a single Ostrovsky equation and we use the method in Appendix B, both of which are similar to schemes in \cite{Khusnutdinova19}. We use a modified pseudo-spectral scheme for the coupled Ostrovsky equations, based upon \cite{Trefethen00}, to allow for larger time steps to be taken, and this is presented in Appendix C. We note that the choice of integration constants to maintain zero mean is conveniently implemented within the pseudo-spectral scheme by setting the coefficient of the zero harmonic to zero.

Unless stated otherwise, we assume that $\Delta x = 0.1$, $\Delta t = 0.01$ and $\Delta T = \varepsilon \Delta t$. In some calculations these parameters may be changed to obtain a higher accuracy result and this will be stated in the figure captions. The domain is taken as $[-L, L]$ in all cases, with $L = 300$ for most calculations with solitary waves and $3 \times L_K$ for cnoidal waves, where $L_K$ is the period of the cnoidal wave.

%%%%%%%%%%%%%%%%%%%%%%%%%%%%%%%%%%%%%%%%%%%%%%%%%%%%%%%%%%%%%%%%%%%%%%%%%%%%%%%%

In all subsequent calculations, we shall refer to the solution of the system \eqref{ueq} - \eqref{weq} as the numerical solution of the exact system and compare this to the weakly-nonlinear solution with an increasing number of terms included. We calculate the solution in the domain $x \in [-40, 40]$ and for $t \in [0, \hat{T}]$ where $\hat{T} = 1/\varepsilon$. 

To construct a right-propagating wave as the initial condition, with no left-propagating wave, we choose our functions $F$ and $V$ appropriately so that $f_1^{-} = F_1$, $f_2^{-} = F_2$ and $f_1^+ = f_2^+ = 0$. Explicitly, for solitary wave initial conditions, we have
\begin{align}
	F_{1}(x) &= A_{1} \sechn{2}{\frac{x}{\Lambda_{1}}} + p, \quad F_{2}(x) = A_{2} \sechn{2}{\frac{x}{\Lambda_{2}}}, \notag \\
	V_{1}(x) &= 2 \frac{A_{1}}{\Lambda{1}} \sechn{2}{\frac{x}{\Lambda_{1}}} \tanh{\lb \frac{x}{\Lambda_{1}} \rb}, \quad V_{2}(x) = 2 c \frac{A_{2}}{\Lambda_{2}} \sechn{2}{\frac{x}{\Lambda_{2}}} \tanh{\lb \frac{x}{\Lambda_{2}} \rb},
	\label{uwIC}
\end{align}
where $p$ is a constant and we have
\begin{equation*}
	A_{1} = 6k_{1}^2, \quad \Lambda_{1} = \frac{\sqrt{2}}{k_{1}}, \quad k_{1} = \frac{1}{\sqrt{6}}, \quad A_{2} = \frac{6ck_{2}^2}{\alpha}, \quad \Lambda_{2} = \frac{\sqrt{2c\beta}}{k_{2}}, \quad k_{2} = \sqrt{\frac{\alpha}{6c}}.
\end{equation*}
We have added a pedestal to the initial condition for $u$ to have distinct non-zero values for $d_{1}$ and $d_{2}$.

To determine the agreement between the numerical solution of the exact system and the weakly-nonlinear solution, we calculate the error as
\begin{equation}
	e^{(i)} = \log \left| 1 - \frac{v_i}{v} \right|, \quad i = 1, 2, 3,
	\label{NumError}
\end{equation}
where $v$ is the numerical solution of the exact system, the weakly-nonlinear solution \eqref{uWNLc1}, \eqref{wWNLc1} with only leading order terms included as $v_{1}$, with terms up to and including $\O{\sqrt{\varepsilon}}$ as $v_{2}$ and with terms up to and including $\O{\varepsilon}$ as $v_{3}$. This will be plotted alongside $v$ for comparison.

\begin{figure}[!htbp]
	\begin{center}
		\includegraphics[width=0.5\textwidth]{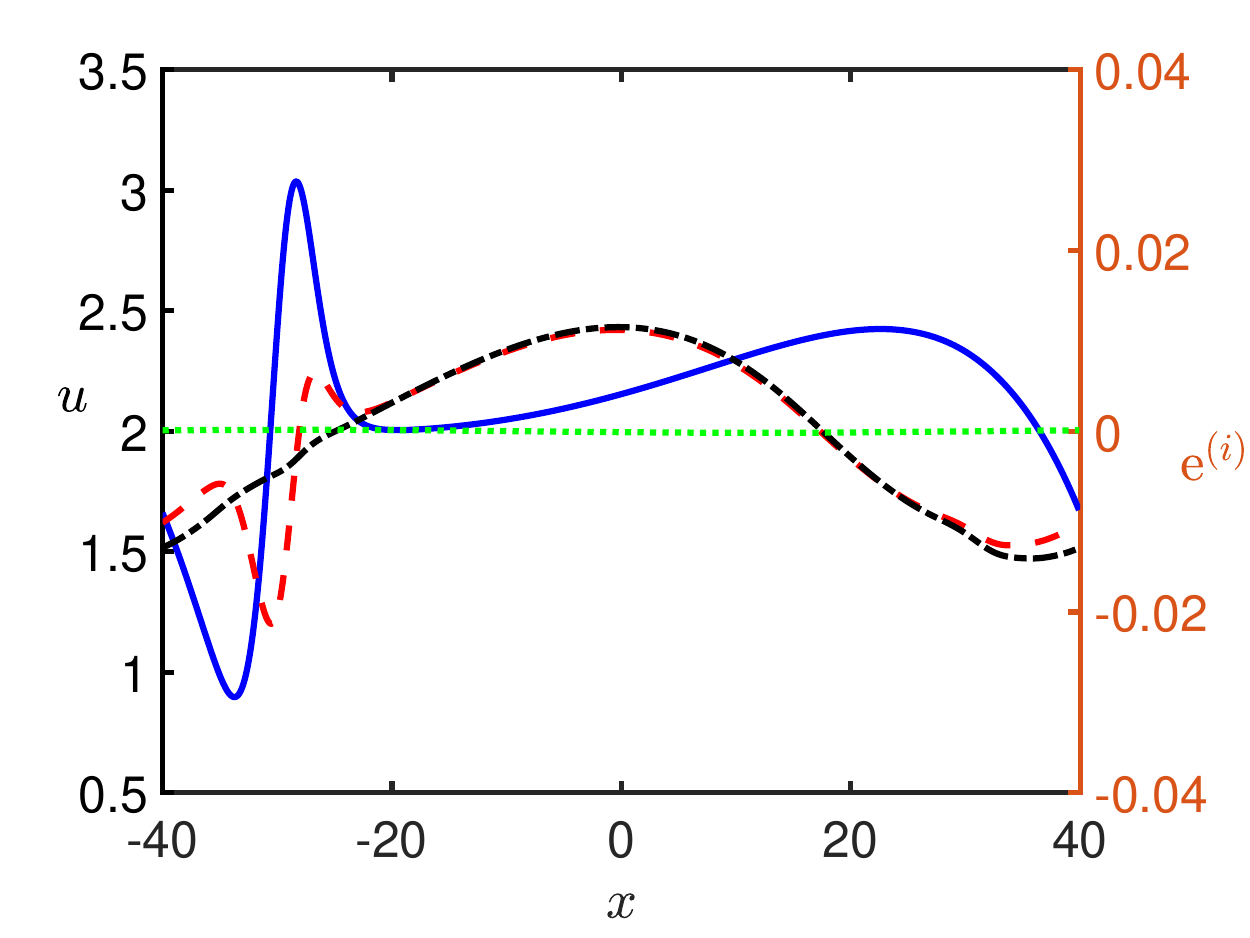}
		\caption{\small A result of the direct numerical simulation (solid, blue) and comparison of the errors between this solution and the weakly-nonlinear solution, including leading order (dashed, red), $\O{\sqrt{\varepsilon}}$ (dash-dot, black) and $\O{\varepsilon}$ (dot, green) corrections, at $t=1/\varepsilon$, for $u$. Parameters are $L=40$, $N=800$, $k = 1/\sqrt{6}$, $\alpha = \beta = c = 2$, $\gamma = 0.5$, $p = 7$, $\varepsilon = 2.38 \times 10^{-3}$, $\Delta t = 0.01$ and $\Delta T = \varepsilon \Delta t$. The solution agrees well to leading order and this agreement is improved with the addition of higher-order corrections.}
		\label{fig:ErrPlotsc1}
	\end{center}
\end{figure}

We choose $\alpha = \beta = c = 2$, $\delta = \gamma = 0.5$, $p = 7$ and present the comparison between the numerical solution of the exact system and weakly-nonlinear solutions, at various orders of $\varepsilon$, for $u$ in Figure \ref{fig:ErrPlotsc1}. A similar result can be observed for $w$.
We see from the error lines that the leading order solution (red, dashed line) is reasonably accurate, with an error of less than 0.02. This is improved with the addition of the $\O{\sqrt{\varepsilon}}$ terms (black, dash-dotted line), reducing the phase shift around the main wave packet. The inclusion of $\O{\varepsilon}$ terms (green, dotted line) reduces the error by an order of magnitude, as expected. We only present one example here for brevity, but the same results have been observed for multiple values of $\delta$, $\gamma$, pedestal height $p$, and for $w$.

To further confirm the validity of the solution, we calculate the maximum absolute error over $x$ as
\begin{equation}
	r_{i} = \max_{-L \leq x \leq L} \lvert v \lb x, t \rb - v_{i} \lb x, t \rb \rvert, \quad i = 1, 2, 3.
	\label{MaxErr}
\end{equation}
This error is calculated at every time step and, to smooth the oscillations in the errors we average $r_{i}$ in the final third of the calculation, denoting this value as $\hat{r}_{i}$. We then use a least-squares power fit to determine how the maximum absolute error varies with the small parameter $\varepsilon$. Therefore we write the errors in the form
\begin{equation}
	\mathrm{exp} \lsq \hat{r}_{i} \rsq = C_{i} \varepsilon^{\alpha_{i}},
	\label{Err}
\end{equation}
and take the logarithm of both sides to form the error plot.

\begin{figure}[!htbp]
\begin{center}
	\includegraphics[width=0.5\textwidth]{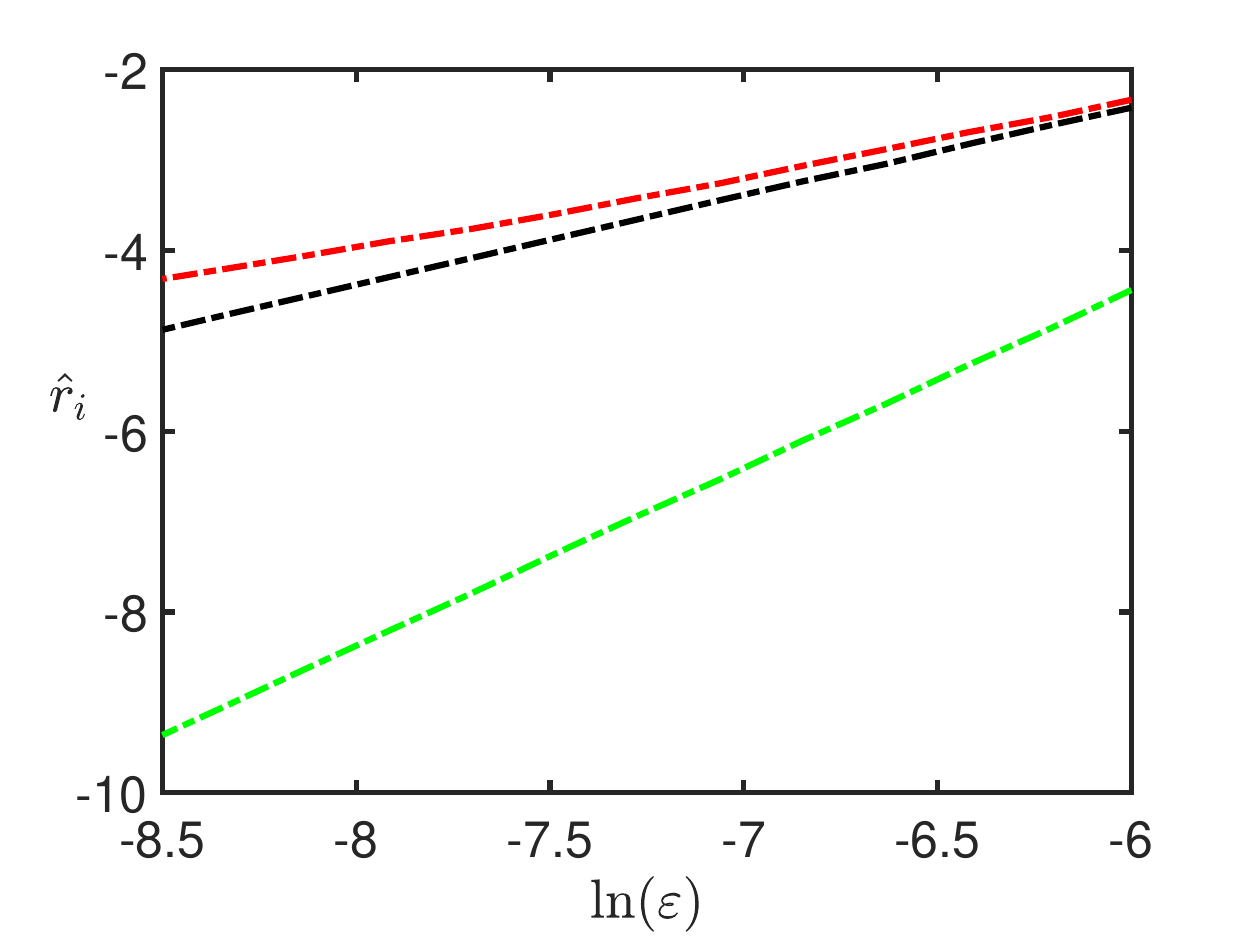}
	\caption{\small A comparison of error curves for varying values of $\varepsilon$, at $t=1/\varepsilon$, for the weakly-nonlinear solution including leading order (upper, red), $\O{\sqrt{\varepsilon}}$ (middle, black) and $\O{\varepsilon}$ (lower, green) corrections for $u$. Parameters are $L=40$, $N=800$, $k = 1/\sqrt{6}$, $\alpha = \beta = c = 2$, $\gamma = 0.5$, $\Delta t = 0.01$ and $\Delta T = \varepsilon \Delta t$. The inclusion of more terms in the expansion increases the accuracy.}
	\label{fig:ErrCompc1}
	\end{center}
\end{figure}

The corresponding errors for the cases considered in Figure \ref{fig:ErrPlotsc1} are plotted in Figure \ref{fig:ErrCompc1}. 
The slope of the curves is approximately 0.77, 0.99 and 1.96, which are slightly higher than the theoretical values for the inclusion of leading order terms and the case with terms up to and including $\O{\varepsilon}$ terms. However this can be explained by the solutions in this case being wave packets, so a phase shift at the level of $\O{\sqrt{\varepsilon}}$ will not have as much as an effect on the errors as it would for radiating solitary waves.\cite{Khusnutdinova19Book} Furthermore, the increase in the slope is consistent with the results seen for previous studies for the Boussinesq-Klein-Gordon equation.\cite{Khusnutdinova19}

% Cases for co- and counter-propagating solitary waves
\section{Counter-Propagating Waves}
\label{sec:Counter}
We now consider various cases for wave interaction, where we will use the constructed weakly-nonlinear solution with the inclusion of terms up to $\O{\sqrt{\varepsilon}}$ to take account of the mass. Two types of initial condition are considered: solitary wave and cnoidal wave. While our focus is on the high-contrast case (when $c - 1 = \O{1}$), it is instructive to compare with the solutions emerging in the low-contrast case (when $c - 1 = \O{\varepsilon}$). The results for the low-contrast case are presented using the solution obtained in \cite{Khusnutdinova19Book}. Then, we discuss the solutions for the high-contrast case, whose weakly-nonlinear solution was constructed in Section \ref{sec:WNL}.

\subsection{Solitary Wave Initial Condition}
\label{sec:CounterSol}
Firstly we take solitary wave initial conditions, as was done in Section \ref{sec:Validity}. For counter-propagating waves, we introduce a second wave, with the same parameters as the first wave, at a phase shift with a sign change in $V$ to result in wave propagation to the left for the second wave. Explicitly we have
\begin{align}
	F_{i}(x) &= A_{i} \sechn{2}{\frac{x + x_0}{\Lambda_{i}}} + A_{i} \sechn{2}{\frac{x + x_1}{\Lambda_{i}}}, \notag \\
	V_{i}(x) &= 2 \frac{A_{i}}{\Lambda{i}} \sechn{2}{\frac{x + x_0}{\Lambda_{i}}} \tanhn{ }{\frac{x + x_0}{\Lambda_{i}}} - 2 \frac{A_{i}}{\Lambda{i}} \sechn{2}{\frac{x + x_1}{\Lambda_{i}}} \tanhn{ }{\frac{x + x_1}{\Lambda_{i}}},	\label{uwICCounter}
\end{align}
where $i = 1$ is the initial condition for $u$ and $i = 2$ is the initial condition for $w$ and $A_i$, $\Lambda_i$ are as defined in \eqref{uwIC}. These can then be used to give initial conditions for the constructed weakly-nonlinear solution via \eqref{f1ICc1} and \eqref{f2ICc1}.

% RSW
\subsubsection{Low-contrast case}
\begin{figure}[!htbp]
	\begin{center}
	\subfigure[Weakly-nonlinear solution to cRB equations at $t = 150$.]{\includegraphics[width=0.45\textwidth]{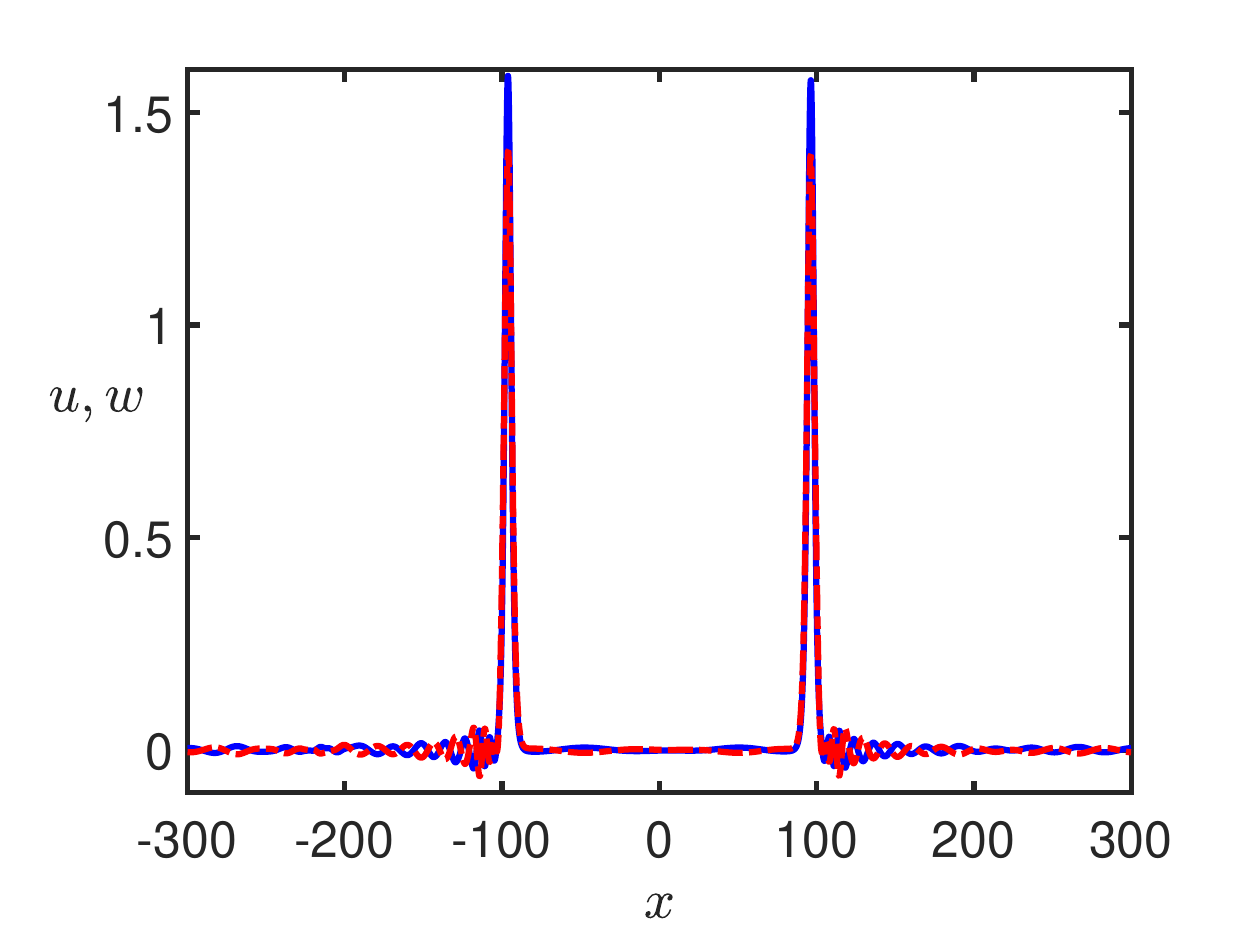}}
	\subfigure[Weakly-nonlinear solution to cRB equations at $t = 500$.]{\includegraphics[width=0.45\textwidth]{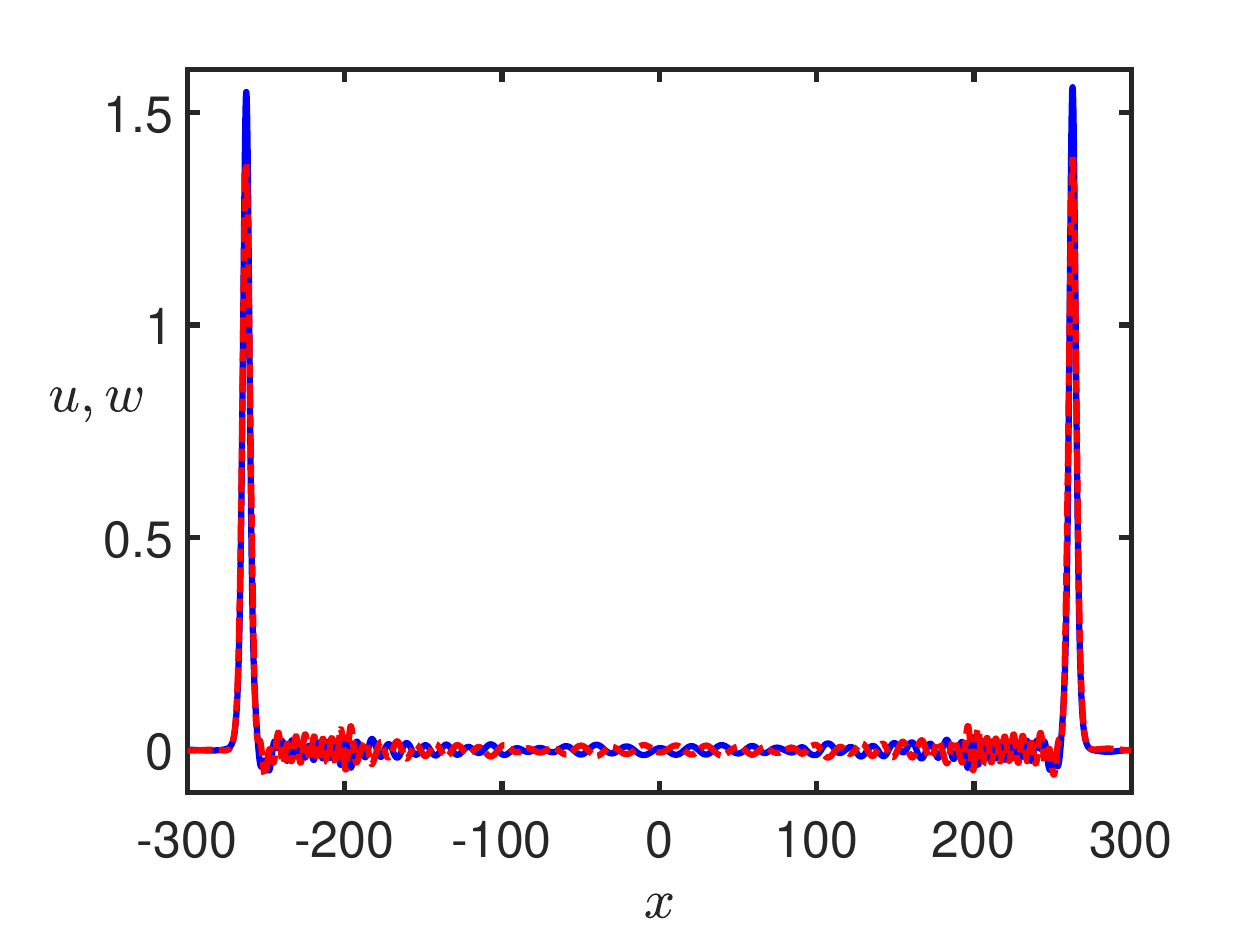}}
	\caption{\small Solution for counter-propagating waves in the cRB equations, using the weakly-nonlinear solution including terms up to $\O{\sqrt{\varepsilon}}$, at (a) $t = 150$ and (b) $t = 500$. Parameters are $L = 300$, $N = 6000$, $\varepsilon = 0.05$, $\alpha = \beta = 1.05$, $c = 1.025$, $\delta = \gamma = 1$, $k_1 = 1/\sqrt{2}$, $k_2 = \sqrt{1.05/2}$, $x_0 = 250$, $x_1 = -250$, $\Delta t = 0.01$ and $\Delta T = \varepsilon \Delta t$. Here, $u/w$ are shown by blue, solid/red, dashed curves, respectively.}
	\label{fig:RSWCounter}
	\end{center}
\end{figure}

We analyse the case when the characteristic speeds in the equations are close, so we have $c - 1 = \O{\varepsilon}$ and the periodic domain is sufficiently large. The solutions in this case are expected to be radiating solitary waves \cite{Khusnutdinova19Book}. We take the initial condition \eqref{uwICCounter}, where the waves will be well separated via the choice of $x_0$ and $x_1$. To determine the behaviour of the solution, we consider it at two times: $t = 150$, before interaction, and $t = 500$, after the interaction. The results are shown in Figure \ref{fig:RSWCounter}.
Firstly we see that the solitons form a solitary wave with a co-propagating radiating tail, as expected from the co-propagating case \cite{Khusnutdinova19Book}. As they interact the solitons appear to emerge with only a small change in their amplitudes, although the amplitude can decay in the evolution of the radiating solitary wave. To determine if their collision is elastic or does indeed introduce a small change in amplitude, or perhaps a phase shift, we compare the solution in the counter-propagating case to a corresponding one when there is only one radiating solitary wave present.

\begin{figure}[!htbp]
	\begin{center}
	\subfigure[Left-propagating radiating solitary wave for $u$ at $t = 500$.]{\includegraphics[width=0.45\textwidth]{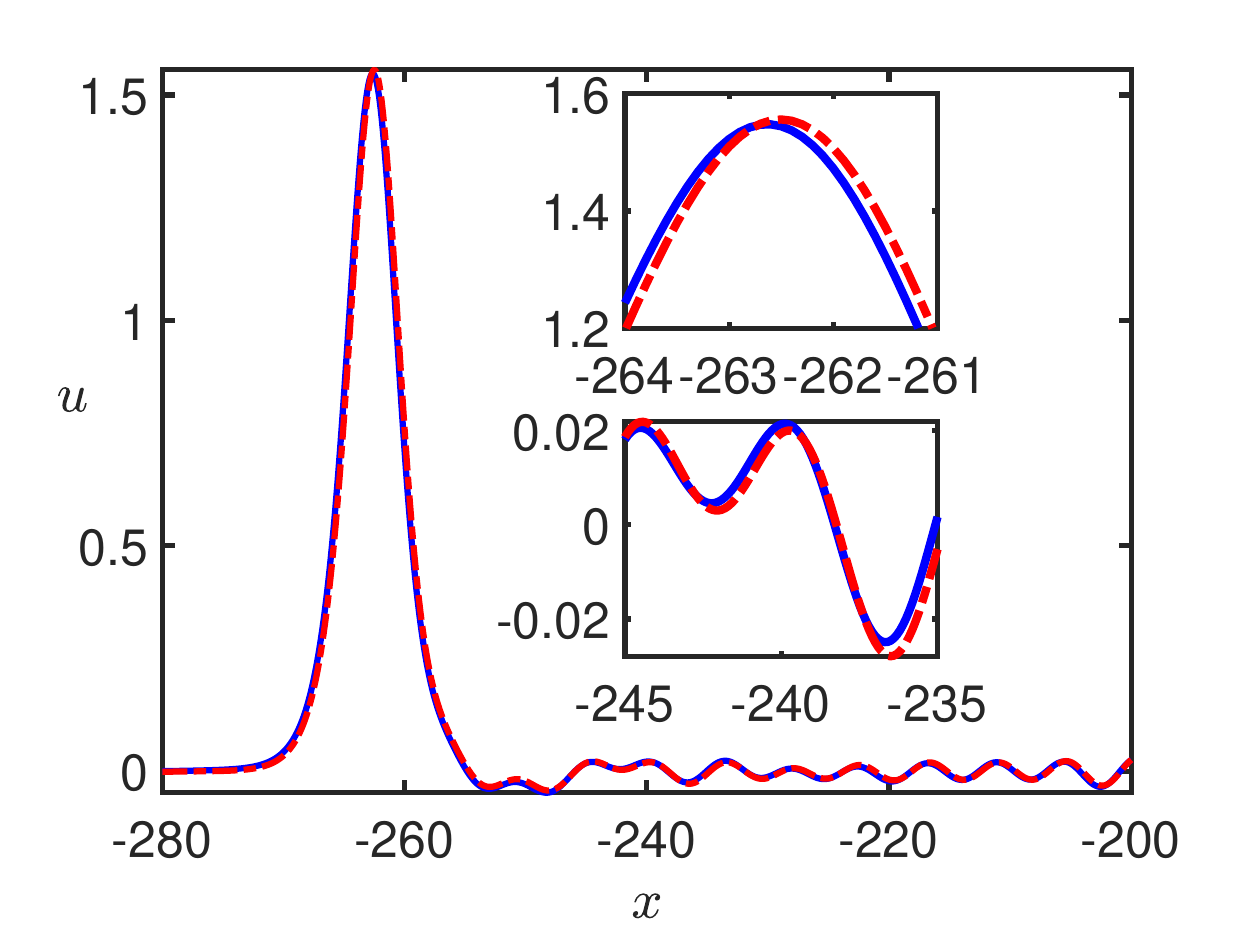}}
	\subfigure[Right-propagating radiating solitary wave for $u$ at $t = 500$.]{\includegraphics[width=0.45\textwidth]{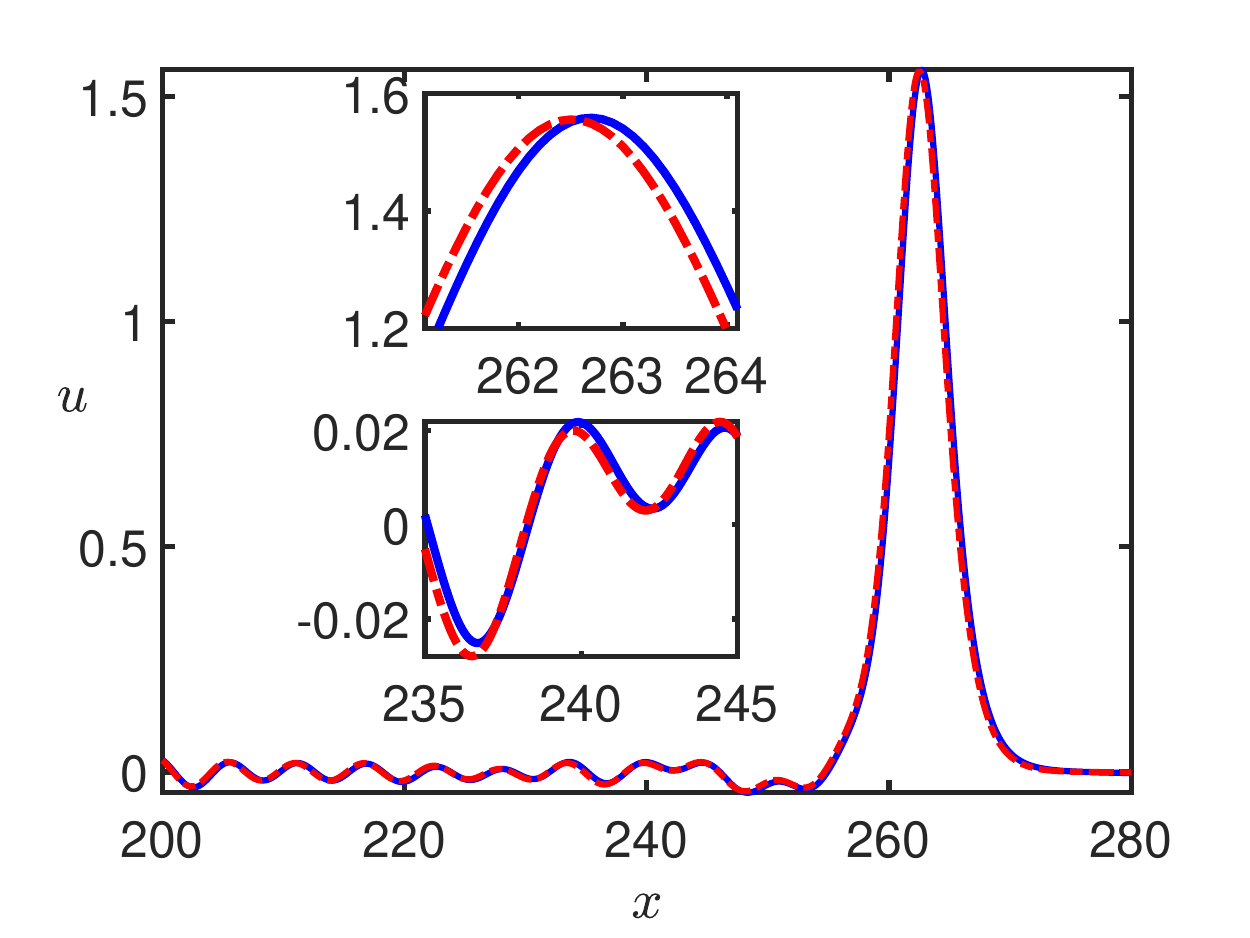}}
	\caption{\small A small phase shift appears as a result of the interaction of radiating solitary waves (blue, solid line) and the case without interaction (red, dashed line), for $u$ at $t = 500$, enhanced for the (a) left-propagating wave, and (b) right-propagating wave. Parameters are $L = 300$, $N = 6000$, $\varepsilon = 0.05$, $\alpha = \beta = 1.05$, $c = 1.025$, $\delta = \gamma = 1$, $k_1 = 1/\sqrt{2}$, $k_2 = \sqrt{1.05/2}$, $x_0 = 250$, $x_1 = -250$, $\Delta t = 0.01$ and $\Delta T = \varepsilon \Delta t$.}
	\label{fig:RSWCounterComp}
	\end{center}
\end{figure}

The comparison between the cases with and without interaction are shown in Figure \ref{fig:RSWCounterComp} for $u$, while a similar result is seen for $w$ and so is omitted for brevity. We can clearly see that a small phase shift occurs once interaction has taken place (blue, solid line) and the amplitude is also slightly reduced on the peak, suggesting that this interaction is not elastic, in contrast to the case with solitary waves. However, the difference is small, and the interaction is only weakly-inelastic.

% Ostrovsky wave packet
\subsubsection{High-contrast case}
\begin{figure}[!htbp]
	\begin{center}
	\subfigure[Solution for $u$ at multiple times.]{\includegraphics[width=0.45\textwidth]{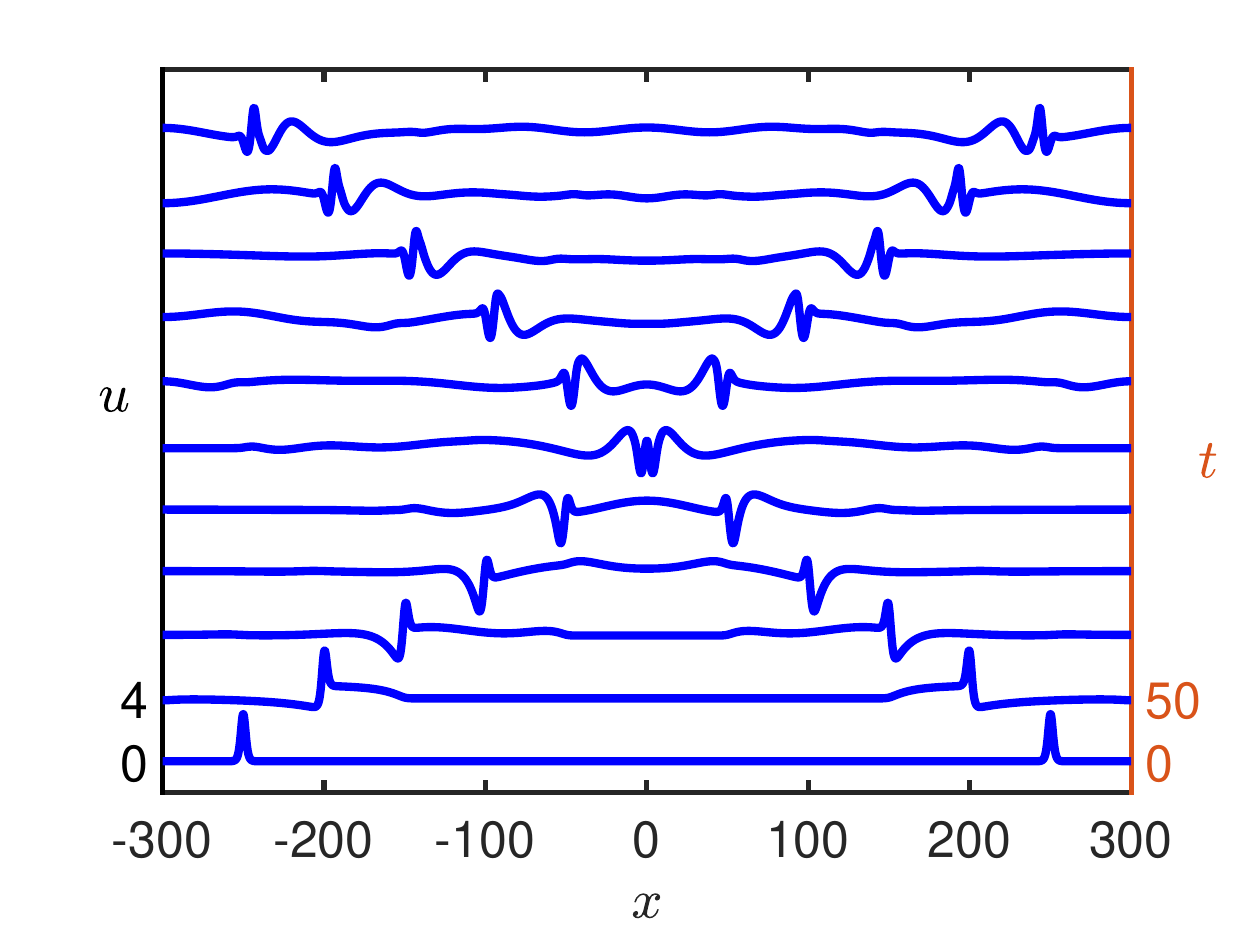}}
	\subfigure[Solution for $w$ at multiple times. The wave moves twice as fast in this case.]{\includegraphics[width=0.45\textwidth]{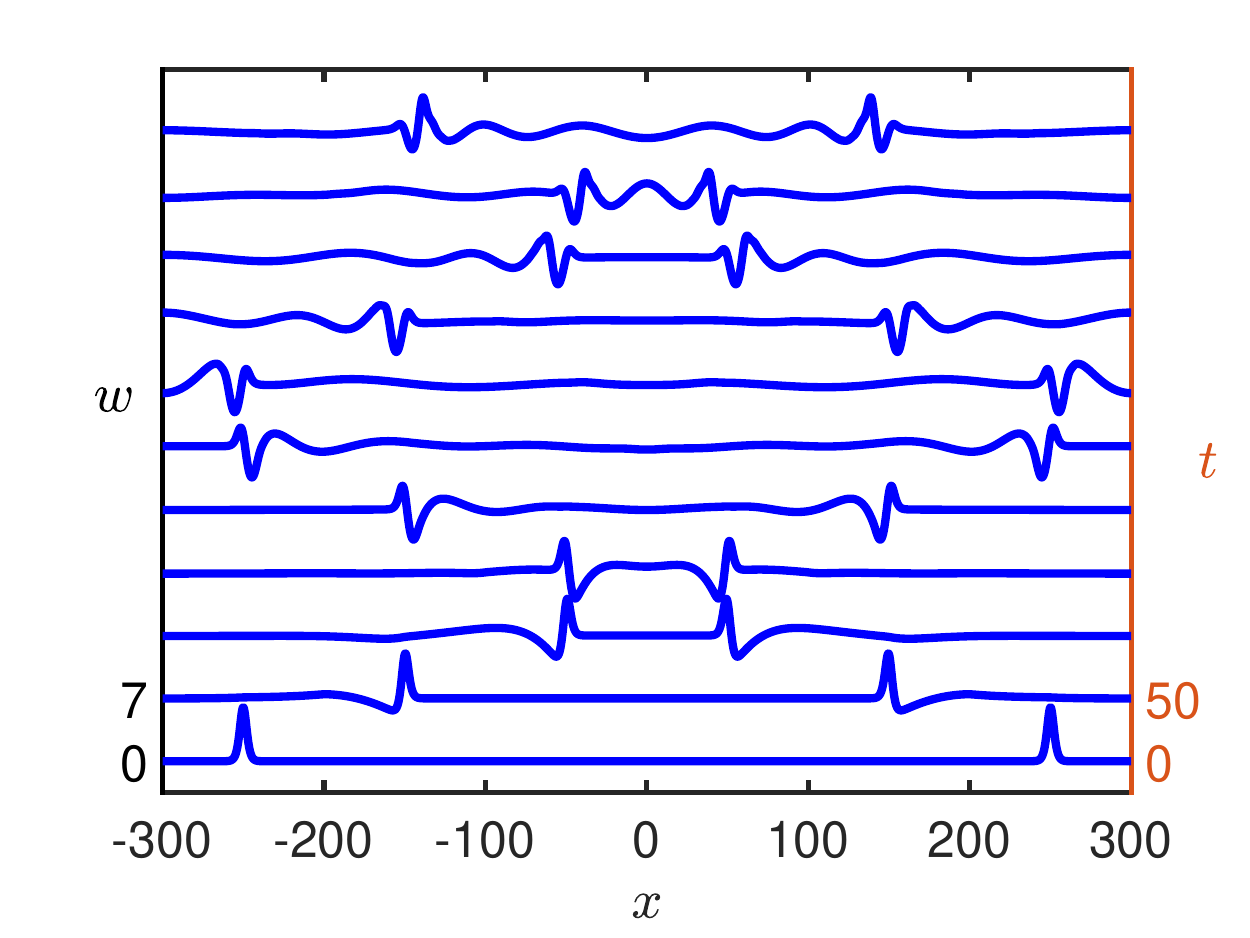}}
	\caption{\small Comparison of counter-propagating wave packets at multiple times, for (a) $u$ and (b) $w$. Parameters are $L = 300$, $N = 6000$, $\varepsilon = 0.0005$, $\alpha = \beta = c = 2$, $\delta = \gamma = 1$, $k_1 = 1/\sqrt{2}$, $k_2 = 1$, $x_0 = 250$, $x_1 = -250$, $\Delta t = 0.01$ and $\Delta T = \varepsilon \Delta t$.}
	\label{fig:OstCounter}
	\end{center}
\end{figure}

Now we consider the case when the characteristic speeds in the equations are distinct, so we have $c - 1 = \O{1}$, and we consider the effect of wave interaction on wave packets. The initial condition is taken as \eqref{uwICCounter} and we analyse the interaction of the generated wave packets. The results are presented in Figure \ref{fig:OstCounter}.
Wave packets are generated in both layers with the packet generated in the lower layer moving faster than the corresponding packet in the upper layer. At the point of interaction the packets move through each other and emerge with the appearance of small changes in their shape and structure, which could be attributed to evolution of the wave packets. 

\begin{figure}[!htbp]
	\begin{center}
	\subfigure[Left-propagating wave packet for $u$ at $t = 500$.]{\includegraphics[width=0.45\textwidth]{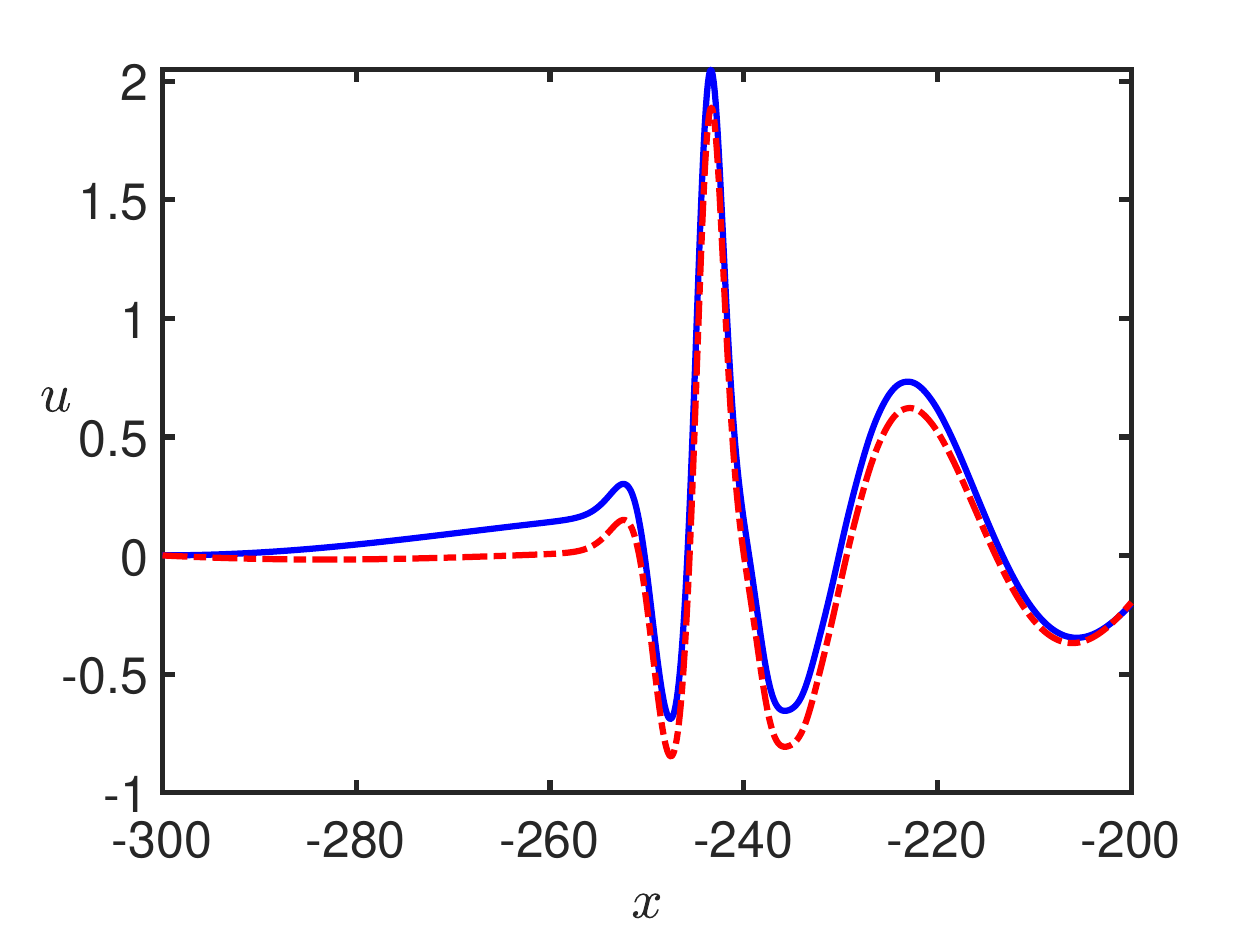}}
	\subfigure[Right-propagating wave packet for $u$ at $t = 500$.]{\includegraphics[width=0.45\textwidth]{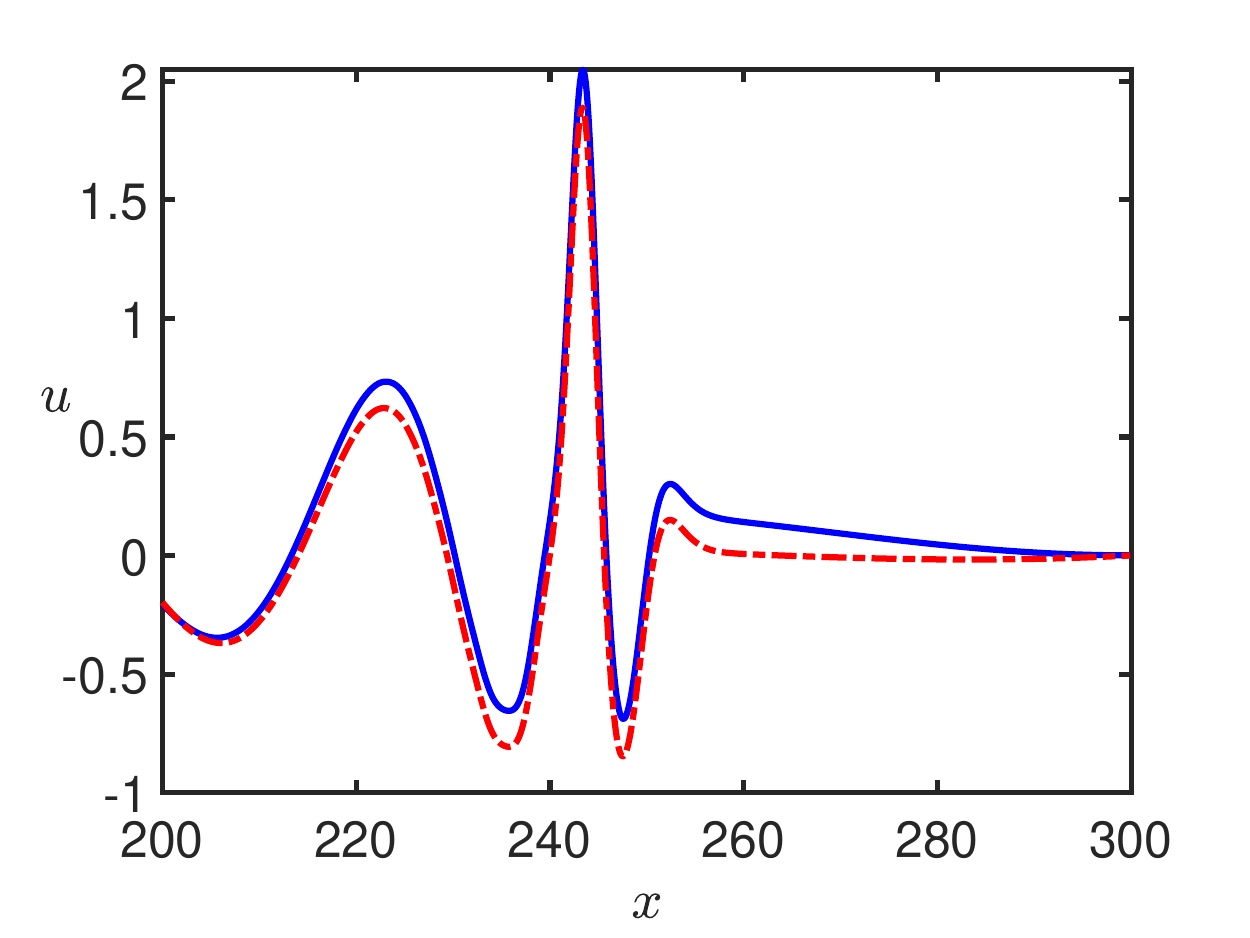}}
	\caption{\small Significant differences appear as a result of the interaction of wave packets (blue, solid line) and the case without interaction (red, dashed line), for $u$ at $t = 500$, enhanced for the (a) left-propagating wave, and (b) right-propagating wave.  Parameters are $L = 300$, $N = 6000$, $\varepsilon = 0.0005$, $\alpha = \beta = c = 2$, $\delta = \gamma = 1$, $k_1 = 1/\sqrt{2}$, $k_2 = 1$, $x_0 = 250$, $x_1 = -250$, $\Delta t = 0.01$ and $\Delta T = \varepsilon \Delta t$.}
	\label{fig:OstCounterComp}
	\end{center}
\end{figure}

As was done in the low-contrast case, we compare the result of the interaction to the corresponding co-propagating case without interaction. The comparison is shown in Figure \ref{fig:OstCounterComp} for $u$, and again we see a similar result for $w$, albeit at a different position due to the different characteristic speeds, so we omit it here. We see that interaction (blue, solid line) leads to several changes, in particular the wave packet is linked to the other wave packet in the case of interaction, leading to a number of differences between the solutions. Therefore, the collision is strongly inelastic in the high-contrast case.

% Study of cnoidal waves
\subsection{Cnoidal Wave Initial Condition}
\label{sec:CounterCnoidal}
A second case of interest is using a cnoidal wave as the initial condition for the cRB equations. We can obtain the cnoidal wave initial condition by considering the uncoupled equations for $u$ and $w$, for example
\begin{equation}
	w_{tt} - c^2 w_{xx} = \varepsilon \lsq \frac{\alpha}{2} \lb w^2 \rb_{xx} + \beta w_{ttxx} \rsq
	\label{wCnoidal}
\end{equation}
for $w$, and the equation for $u$ can be found by setting $c = \alpha = \beta = 1$.
The leading-order weakly-nonlinear solution to this equation takes the form of the KdV equation 
\begin{equation}
	2c f_{T} + \alpha f f_{\xi} + \beta c^2 f_{\xi \xi \xi} = 0,
	\label{CnoidalKdV}
\end{equation}
where we have assumed that $\xi = x - ct$ and $T = \varepsilon t$. This can be thought of as a truncated Ostrovsky equation derived in Section \ref{sec:WNL}, under the assumption that the initial condition has zero mean. The cnoidal wave solution to this equation can be obtained as
\begin{equation}
	f = -\frac{6 \beta c^2}{\alpha} \lb f_2 - (f_2 - f_3) \mathrm{cn}^2 \lsq \lb \xi + \nu T \rb \sqrt{\frac{f_1 - f_3}{2}} | m \rsq  \rb,
	\label{CnoidalIC}
\end{equation}
where
\begin{equation*}
	\nu = (f_1 + f_2 + f_3) \beta c, \quad m = \frac{f_2 - f_3}{f_1 - f_3}.
\end{equation*}
The solution is parametrised by the three constants $f_3 < f_2 < f_1$ with $0 < m < 1$. The wave length of the cnoidal wave is given by
\begin{equation}
	L_K = 2 K(m) \sqrt{\frac{2}{f_1 - f_3}},
	\label{CnoidalWL}
\end{equation}
where 
\begin{equation*}
	K(m) = \int_0^{\pi/2} \frac{d \theta}{\sqrt{1 - m\  \sin^2 \theta}}
\end{equation*}
is the complete elliptic integral of the first kind. This wave is an exact solution to the KdV equation and a natural initial condition for the cRB equations, in the same way as the solitary wave solution in Section \ref{sec:CounterSol} also was a natural initial condition since the Ostrovsky equation is an extension of the KdV equation.

To reduce the number of parameters in the solution, we will assume that $f_2 = 0$, which means that the cnoidal wave is moving on a zero background. This will reduce the magnitude of the mean term, but it will not be zero, which is consistent with the approach for the solitary wave solution. To simplify notation we introduce the constant $\theta = \sqrt{\frac{f_1 - f_3}{2}}$ and variable $\tilde{x} = x + x_0$. 

To construct a right-propagating cnoidal wave as the initial condition, we  choose our functions $F$ and $V$ appropriately, as was done for the solitary wave initial conditions. Explicitly for the co-propagating wave we have the initial condition
\begin{align}
	F_{1}(x) &= -6 f_3\: \mathrm{cn}^2 \lsq \theta \tilde{x} | m \rsq, \quad F_{2}(x) = -\frac{6 \beta c^2}{\alpha} f_3\: \mathrm{cn}^2 \lsq \theta \tilde{x} | m \rsq, \notag \\
	V_{1}(x) &= -12 f_3 \theta\: \mathrm{cn} \lsq \theta \tilde{x} | m \rsq \mathrm{sn} \lsq \theta \tilde{x} | m \rsq \mathrm{dn} \lsq \theta \tilde{x} | m \rsq, \quad	V_{2}(x) = -\frac{12 \beta c^2}{\alpha} f_3 \theta\: \mathrm{cn} \lsq \theta \tilde{x} | m \rsq \mathrm{sn} \lsq \theta \tilde{x} | m \rsq \mathrm{dn} \lsq \theta \tilde{x} | m \rsq.
	\label{uwICCnoidalCo}
\end{align}

\subsubsection{Low-contrast case}
\begin{figure}[!htbp]
	\begin{center}
	\subfigure[Solution at $t = 0$.]{\includegraphics[width=0.45\textwidth]{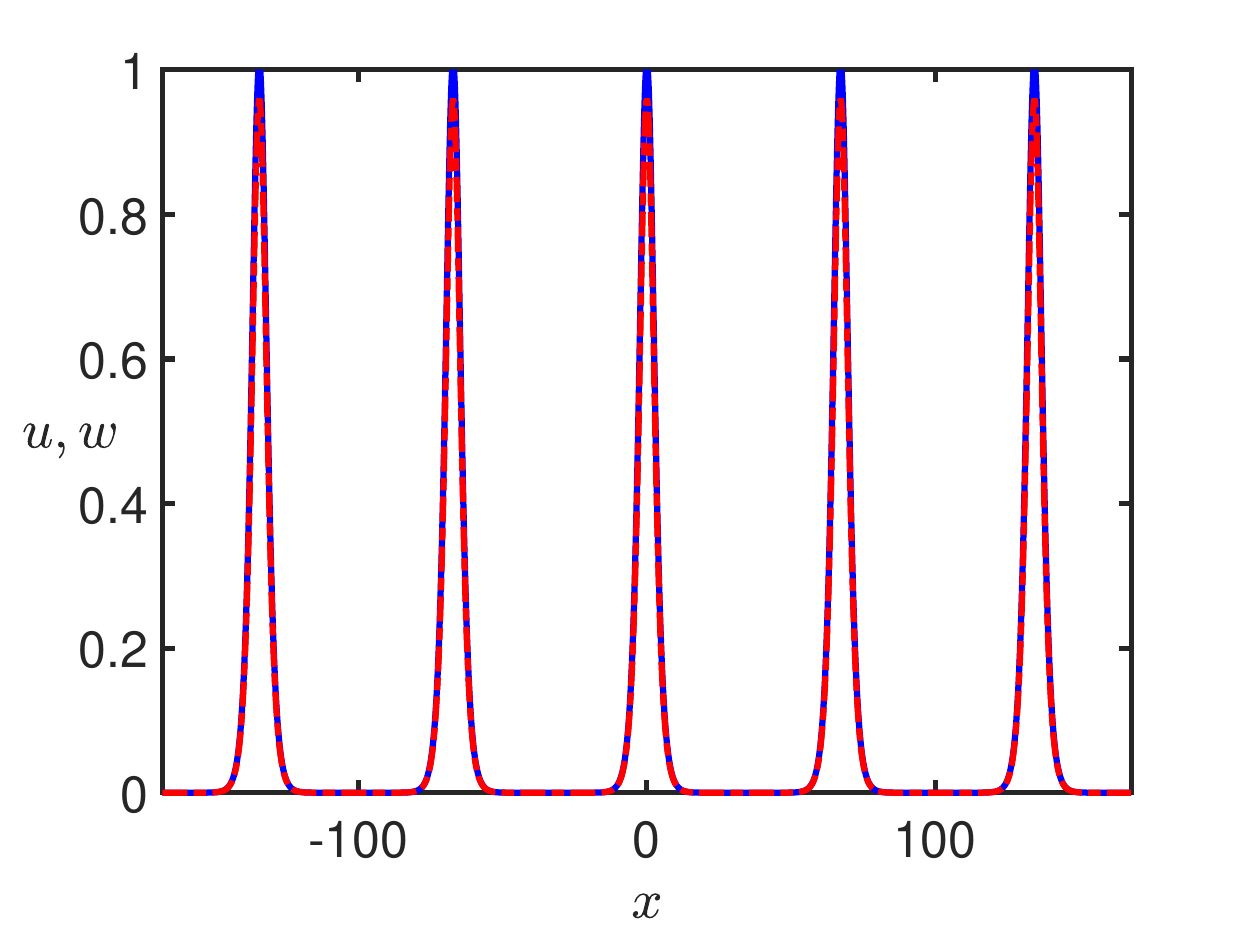}}
	\subfigure[Solution at $t = 400$.]{\includegraphics[width=0.45\textwidth]{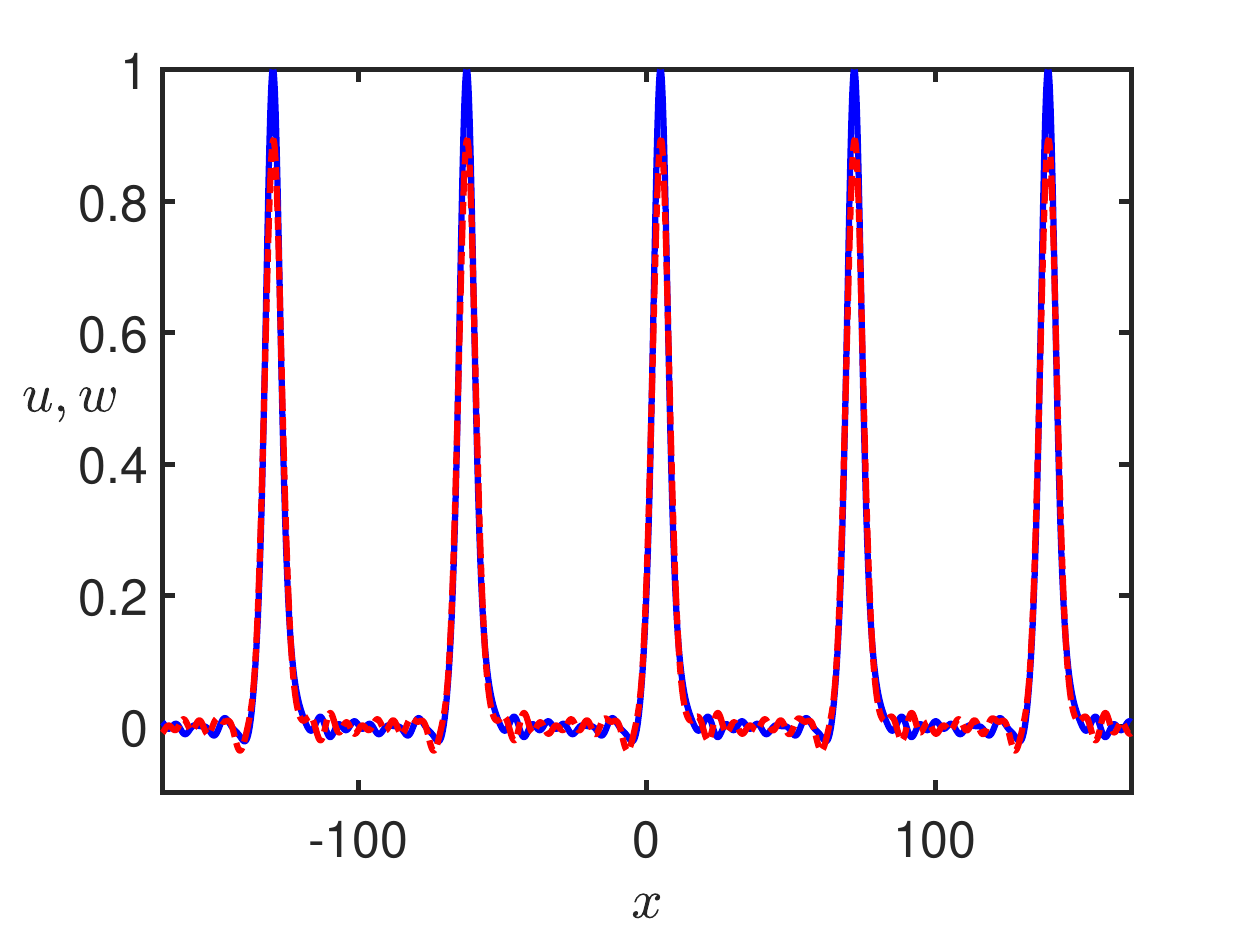}}
	\caption{\small Generation of co-propagating radiating solitary waves from a cnoidal wave initial condition, for $u$ (blue, solid line) and $w$ (red, dashed line), at (a) $t = 0$, and (b) $t = 400$. Parameters are $L \approx 168.0$, $N = 3360$, $\varepsilon = 0.05$, $f_1 = 1 \times 10^{-8}$, $f_2 = 0$, $f_3 = -\frac{1}{6}$, $\alpha = \beta = 1$, $c = 1.025$, $\delta = \gamma = 1$, $x_0 = 0$, $\Delta t = 0.01$ and $\Delta T = \varepsilon \Delta t$.}
	\label{fig:CnoidalCo}
	\end{center}
\end{figure}

Firstly we consider the low-contrast case. Our earlier results would suggest that each peak in the cnoidal wave will evolve into a radiating solitary wave, whose tail will eventually interact with the preceding wave. We take the initial condition for co-propagating waves, given in \eqref{uwICCnoidalCo}, and the results are presented in Figure \ref{fig:CnoidalCo}.

As the troughs of the cnoidal wave are long, we can see the generation of a co-propagating radiating tail forming behind the cnoidal wave peaks, which then begins to interact with the preceding peak. The peaks of the wave appear to survive this interaction and maintain their shape, with a small reduction in amplitude that could be attributed to the generation of the radiating tail. The waves appear to form a larger wave structure across the domain.

\begin{figure}[!htbp]
	\begin{center}
	\subfigure[Solution at $t = 800$ for $u$.]{\includegraphics[width=0.45\textwidth]{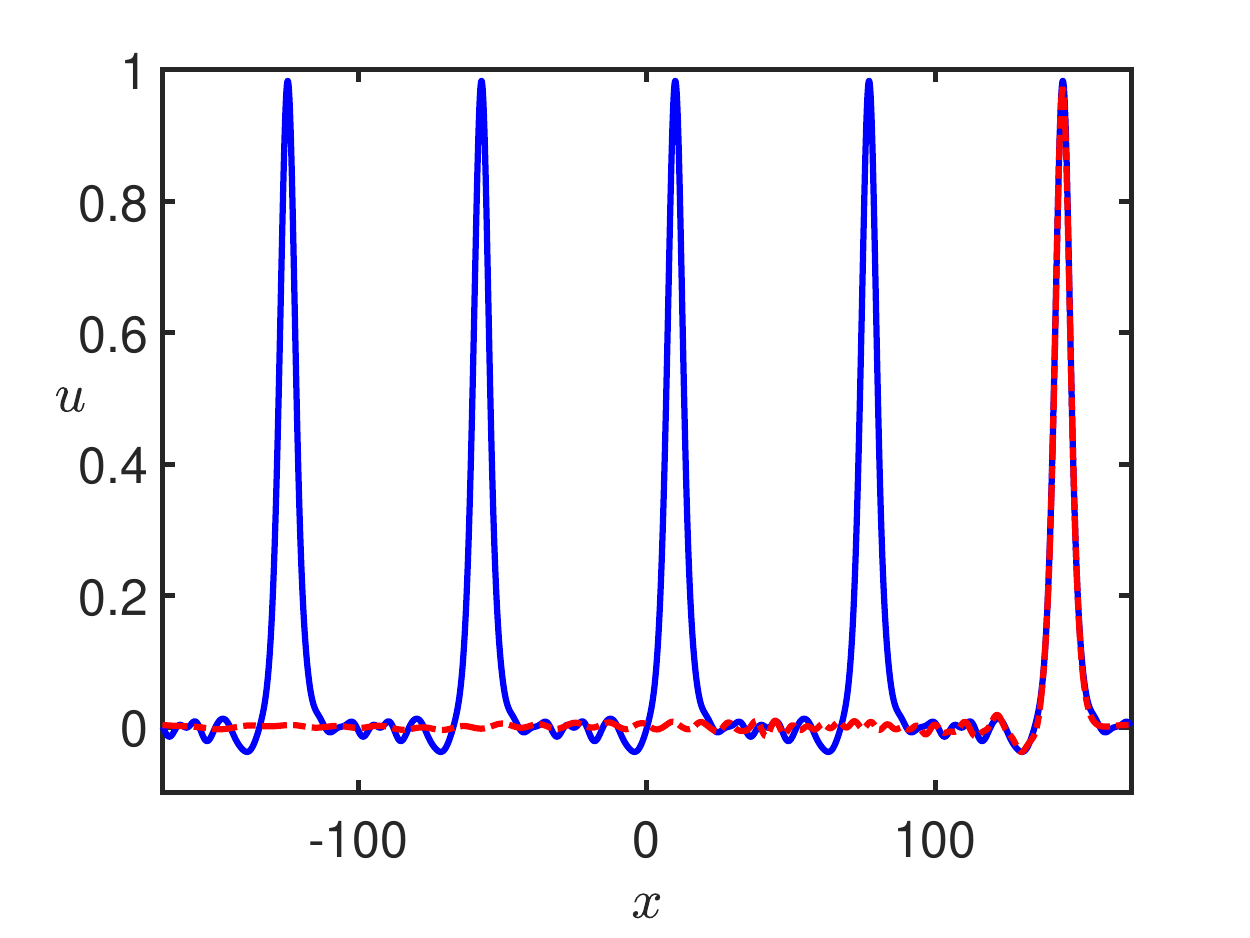}}
	\subfigure[Amplitude of case with (dashed line) and without (solid line) interaction.]{\includegraphics[width=0.45\textwidth]{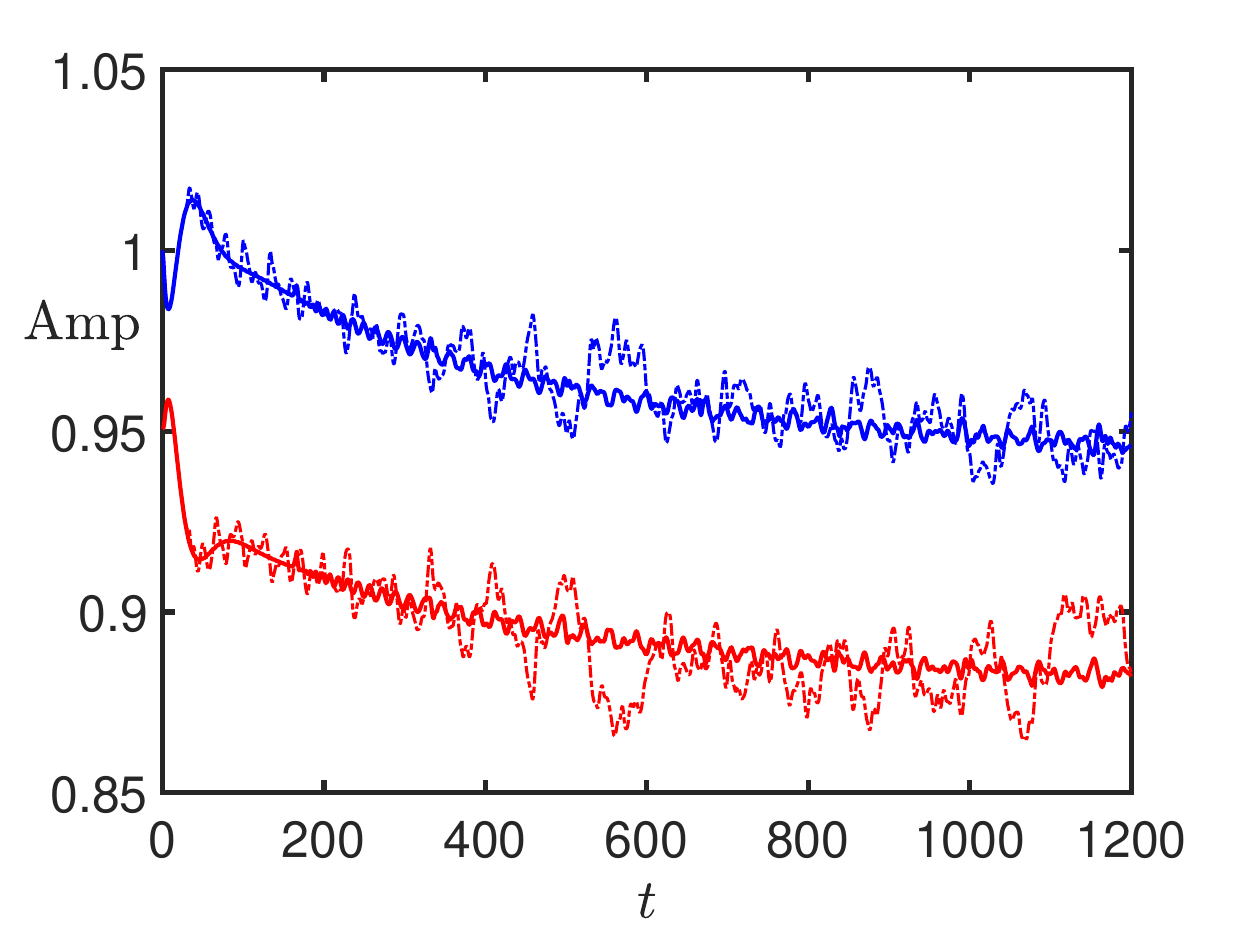}}
	\caption{\small Comparison of effect of peak interaction with radiating tails. In (a) we have a cnoidal wave with multiple peaks (blue, solid line) and a cnoidal wave with a single peak (red, dashed line). In (b) the peak amplitude in the single peak case (solid line) and multiple peak case (dashed line), for $u$ (blue) and $w$ (red). Parameters are $L \approx 168.0$, $N = 3360$, $\varepsilon = 0.05$, $f_1 = 1 \times 10^{-8}$, $f_2 = 0$, $f_3 = -\frac{1}{6}$, $\alpha = \beta = 1$, $c = 1.025$, $\delta = \gamma = 1$, $x_0 = 0$, $\Delta t = 0.01$ and $\Delta T = \varepsilon \Delta t$.}
	\label{fig:CnoidalCoComp}
	\end{center}
\end{figure}

We can consider the effect of interaction of the peaks with the radiating tails by running a simulation for the same initial condition, with five peaks, shown in Figure \ref{fig:CnoidalCo}, and compare this with a truncated initial condition consisting of only a single peak, but with the same domain size as our previous simulation. This is shown in Figure \ref{fig:CnoidalCoComp}.
From Figure \ref{fig:CnoidalCoComp}(a) we can see that there is not a significant phase shift in the leading peak, however we see several phase shifts and amplitude changes have occurred with the tail, caused by the interaction between the tails and the peaks. In Figure \ref{fig:CnoidalCoComp}(b) we can see that the amplitude of the peak in the case without interaction follows a fairly steady path, while in the case with interaction it has large changes caused by interactions with the radiated tail of the preceding peak. This is qualitatively similar to the effects of a large soliton tunelling through a soliton gas or other oscillatory wave structure  (see, for example, \cite{Sande21, Girotti22} and references therein), but in this case there are no significant phase-shifts, and the main effect is seen in the amplitude variations.

\subsubsection{High-contrast case}
\begin{figure}[!htbp]
	\begin{center}
	\subfigure[Solution at $t = 0$.]{\includegraphics[width=0.45\textwidth]{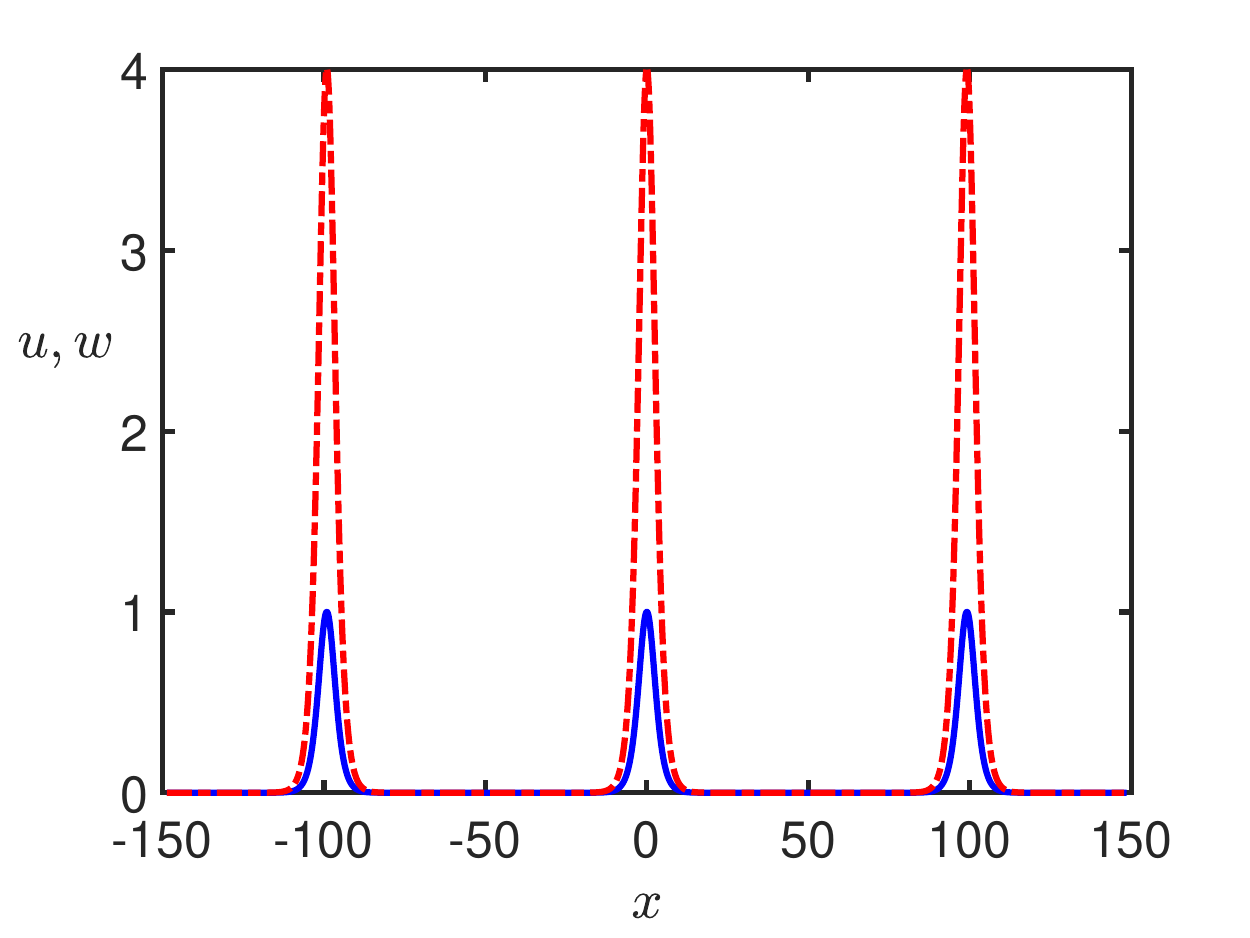}}
	\subfigure[Solution at $t = 500$.]{\includegraphics[width=0.45\textwidth]{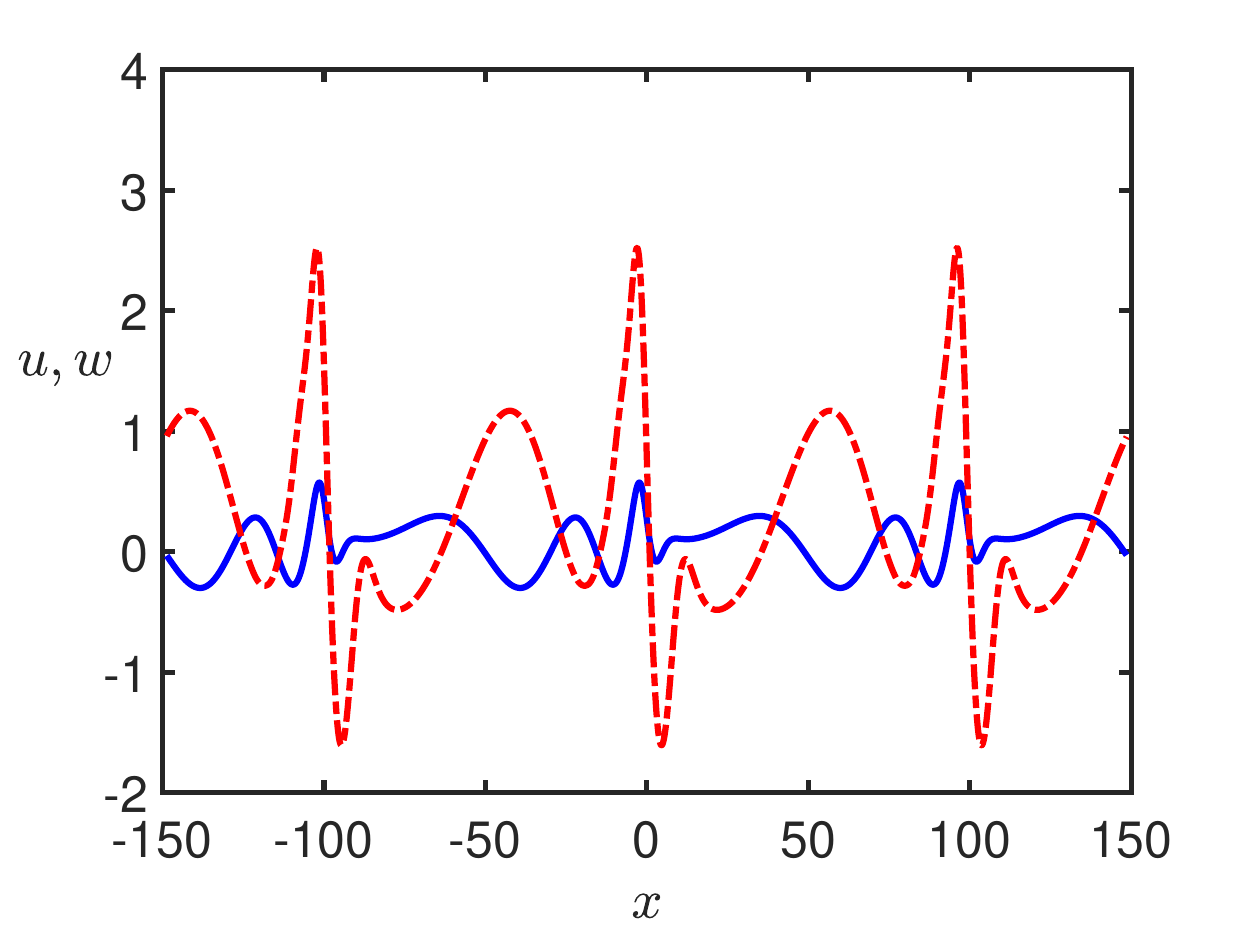}}
	\caption{\small Generation of wave packets from a cnoidal wave initial condition, for $u$ (blue, solid line) and $w$ (red, dashed line), at (a) $t = 0$, and (b) $t = 500$. Parameters are $L \approx 148.6$, $N = 2974$, $\varepsilon = 0.005$, $f_1 = 1 \times 10^{-12}$, $f_2 = 0$, $f_3 = -\frac{1}{6}$, $\alpha = \beta = 2$, $c = 2$, $\delta = \gamma = 1$, $x_0 = 0$, $\Delta t = 0.01$ and $\Delta T = \varepsilon \Delta t$.}
	\label{fig:OstCnoidalCo}
	\end{center}
\end{figure}

Following on from Section \ref{sec:CounterSol}, we now examine the case when the characteristic speeds in the equations are distinct, so we have $c - 1 = \O{1}$. We take the initial condition \eqref{uwICCnoidalCo} and the results are presented in Figure \ref{fig:OstCnoidalCo}.
A wave packet is formed by each peak of the cnoidal wave, which were chosen to be approximately 100 units apart so that the peaks of the cnoidal wave are distinct and each wave packet can be seen clearly. Note that, although the solutions for $u$ and $w$ overlay each other, the packets formed in the lower layer move faster than the upper layer and their overlap is due to the choice of time. The wave packets are connecting to each other in the intervening space between the main wave packets.

\begin{figure}[!htbp]
	\begin{center}
	\includegraphics[width=0.5\textwidth]{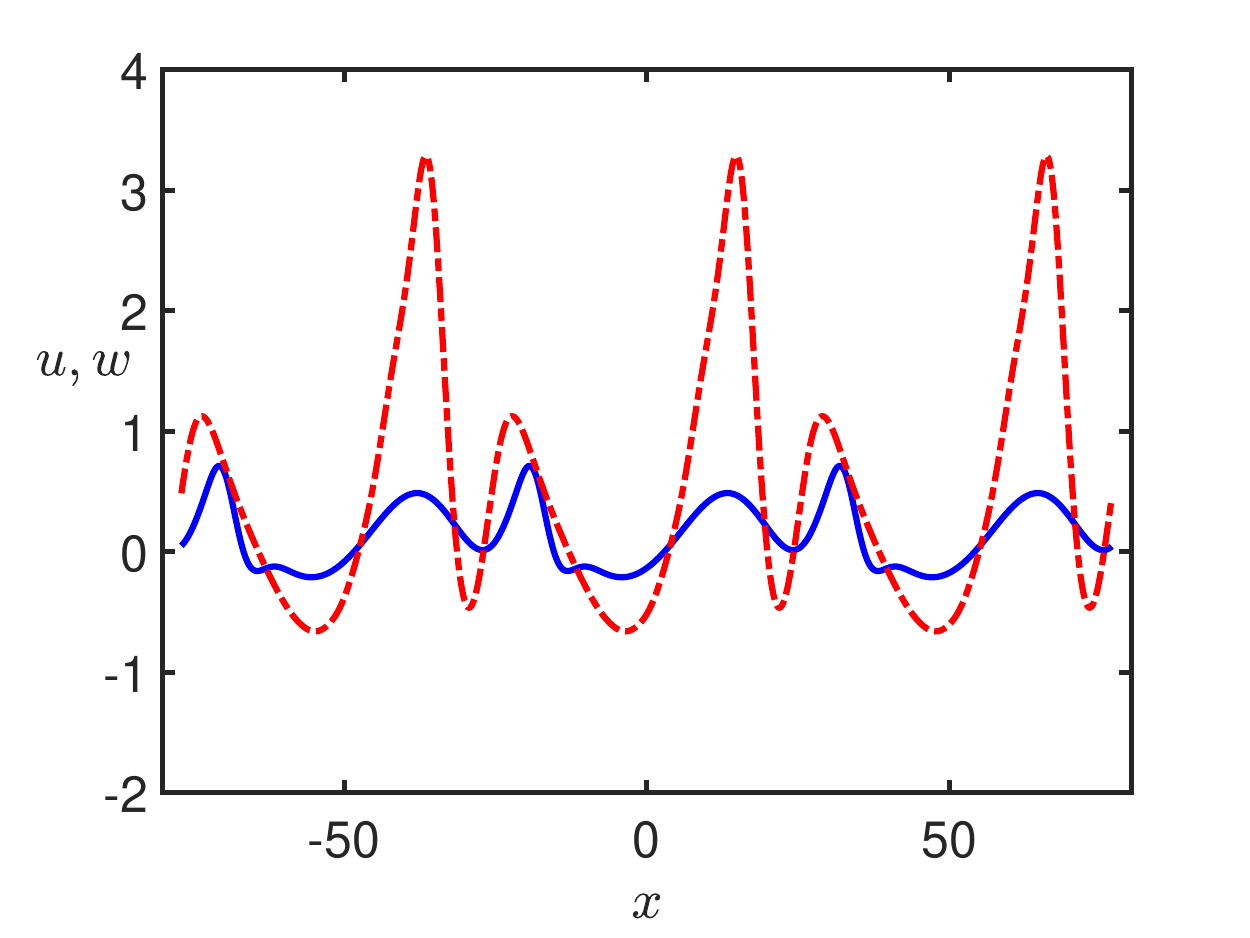}
	\caption{\small Generation of wave packets from a cnoidal wave initial condition, for $u$ (blue, solid line) and $w$ (red, dashed line), at $t = 500$. Parameters are $L \approx 76.9$, $N = 1538$, $\varepsilon = 0.005$, $f_1 = 1 \times 10^{-6}$, $f_2 = 0$, $f_3 = -\frac{1}{6}$, $\alpha = \beta = 2$, $c = 2$, $\delta = \gamma = 1$, $x_0 = 0$, $\Delta t = 0.01$ and $\Delta T = \varepsilon \Delta t$.}
	\label{fig:OstCnoidalCo2}
	\end{center}
\end{figure}

To further explore the evolution of the cnoidal wave initial condition into a series of wave packets, we choose an initial condition when the peaks of the cnoidal wave are approximately 50 units apart. This is plotted in Figure \ref{fig:OstCnoidalCo2}. We can see that wave packets are again generated, however in this case the qualitative structure is different and the connection between the wave packets is shorter.

\section{Conservation laws}
\label{sec:Laws}
The system \eqref{ueq} - \eqref{weq} is related to the system
\begin{align}
	U_{tt} - U_{xx} &= \varepsilon \lsq U_x U_{xx} + U_{ttxx} - \delta \lb U - W \rb \rsq,  \label{Ueq}\\
	W_{tt} - c^2 W_{xx} &= \varepsilon \lsq \alpha W_x W_{xx} + \beta W_{ttxx} + \gamma \lb U - W \rb \rsq, \label{Weq}
\end{align}
\begin{equation*}
	U(x, t) = \int_{-L}^{x} u(\sigma, t) \dd{\sigma}, \quad W(x, t) = \int_{-L}^{x} w(\sigma, t) \dd{\sigma}.
\end{equation*}
It naturally follows that $U(-L) = W(-L) = 0$, and also that $U_t(-L) = W_t(-L) = 0$, which will be used in what follows. The system \eqref{Ueq} - \eqref{Weq} is Lagrangian, and it has three local conservation laws for the mass, energy and momentum.\cite{Khusnutdinova09}
Using these known conservation laws, we can find the conservation laws for our system \eqref{ueq} - \eqref{weq} in the form of one local (mass) and two non-local (energy and momentum) conservation laws. Indeed, differentiating the conservation law for mass from \cite{Khusnutdinova09} with respect to $x$ and rewriting the energy and momentum laws in terms of $u$, $w$ instead of $U_x$, $W_x$ we obtain
	\begin{align}
		{\lb u_t + \frac{\delta}{\gamma} w_t \rb}_t - {\lb u_x + \frac{\delta c^2}{\gamma} w_x + \varepsilon u u_x + \frac{\varepsilon \alpha \delta}{\gamma} w w_x + u_{ttx} + \frac{\varepsilon \beta \delta}{\gamma} w_{ttx} \rb}_x = 0, \label{cl1} \\
		\frac{1}{2} \left \{ U_t^2 + \frac{\delta}{\gamma} W_t^2 + u^2 + \frac{\delta c^2}{\gamma} w^2 + \frac{\varepsilon}{3} \lb u^3 + \frac{\alpha \delta}{\gamma} w^3 \rb + \varepsilon u_{t}^2 + \frac{\varepsilon \beta \delta}{\gamma} w_{t}^2 + \varepsilon \delta (U - W)^2 \right \}_t \nonumber \\
		- \lb U_t u + \frac{\delta c^2}{\gamma} W_t w + \frac{\varepsilon}{2} U_t u^2 + \frac{\varepsilon \alpha \delta}{2 \gamma} W_t w^2 + \varepsilon U_t u_{tt} + \frac{\varepsilon \beta \delta}{\gamma} W_t w_{tt} \rb_x = 0, \label{cl2} \\
		\lb U_t u + \frac{\delta}{\gamma} W_t w + \varepsilon u_{t} u_{x} + \frac{\varepsilon \beta \delta}{\gamma} w_{t} w_{x} \rb_t - \left \{\varepsilon u u_{tt} + \frac{\varepsilon \beta \delta}{\gamma} w w_{tt} \right . \nonumber \\
		+ \left . \frac{1}{2} \left [ U_t^2 + \frac{\delta}{\gamma} W_t^2 + u^2 + \frac{\delta c^2}{\gamma} w^2 + \frac{2 \varepsilon}{3} \lb u^3 + \frac{\alpha \delta}{\gamma} w^3 \rb + \varepsilon u_{t}^2 + \frac{\varepsilon \beta \delta}{\gamma} w_{t}^2 - \varepsilon \delta (U - W)^2 \right ]  \right \}_x = 0. \label{cl3}
	\end{align}
Integrating these conservation laws with respect to $x$ from $-L$ to $L$, using the periodicity of $u$ and $w$ on $[-L, L]$, we obtain one conserved quantity (mass)
\begin{equation}
\diff{ }{t} \lb \langle u \rangle + \frac{\delta}{\gamma} \langle w \rangle \rb = 0, 
\label{mass}
\end{equation}
and two non-local conservation laws (energy and momentum, respectively)
	\begin{align}
		\diff{E}{t} = \frac{1}{2} \diff{ }{t} &\int_{-L}^L \left \{ U_t^2 + \frac{\delta}{\gamma} W_t^2 + u^2 + \frac{\delta c^2}{\gamma} w^2 + \frac{\varepsilon}{3} \lb u^3 + \frac{\alpha \delta}{\gamma} w^3 \rb + \varepsilon u_t^2 + \frac{\varepsilon \beta \delta}{\gamma} w_t^2 + \varepsilon \delta (U - W)^2 \right \} \dd{x} \nonumber \\
		&= \left \{ U_t(L) \lsq u(L) + \frac{\varepsilon}{2} u^2(L) + \varepsilon u_{tt}(L) \rsq + W_t(L) \lsq \frac{\delta c^2}{\gamma} w(L) + \frac{\varepsilon \alpha \delta}{2 \gamma} w^2(L) + \frac{\varepsilon \beta \delta}{\gamma} w_{tt}(L) \rsq \right \}, \label{energy} \\
		\diff{M}{t} = \diff{ }{t} &\int_{-L}^L  \lb U_t u + \frac{\delta}{\gamma} W_t w + \varepsilon u_t u_x + \frac{\varepsilon \beta \delta}{\gamma} w_t w_x \rb \dd{x} = \frac{1}{2} \lsq U_t^2(L) + \frac{\delta}{\gamma} W_t^2(L) - \varepsilon \delta (U(L) - W(L))^2 \rsq. \label{momentum}
	\end{align}
Note that if $\langle u \rangle= \langle w \rangle = 0$ (i.e. $d_1 = d_2 = 0$), then the non-local conservation laws (\ref{energy}) and (\ref{momentum})  yield the usual (local) conservation of energy and momentum. More generally, the non-local conservation laws are applicable when the waves propagate on the background of some initially pre-strained basic state characterised by non-zero mass of $u$ and $w$.

\begin{figure}[!htbp]
	\begin{center}
	\subfigure[Energy for cRB equations.]{\includegraphics[width=0.45\textwidth]{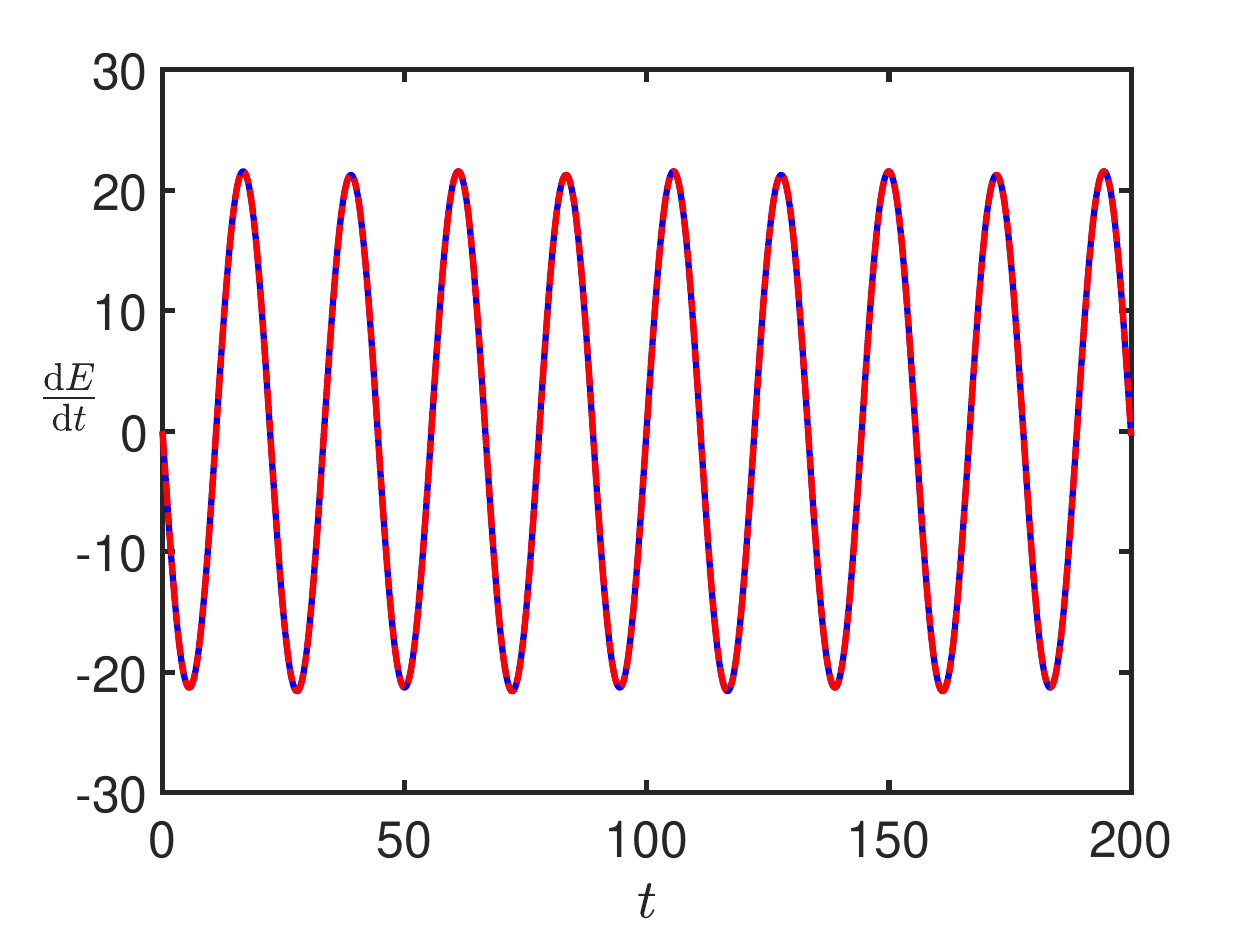}}
	\subfigure[Momentum for cRB equations.]{\includegraphics[width=0.45\textwidth]{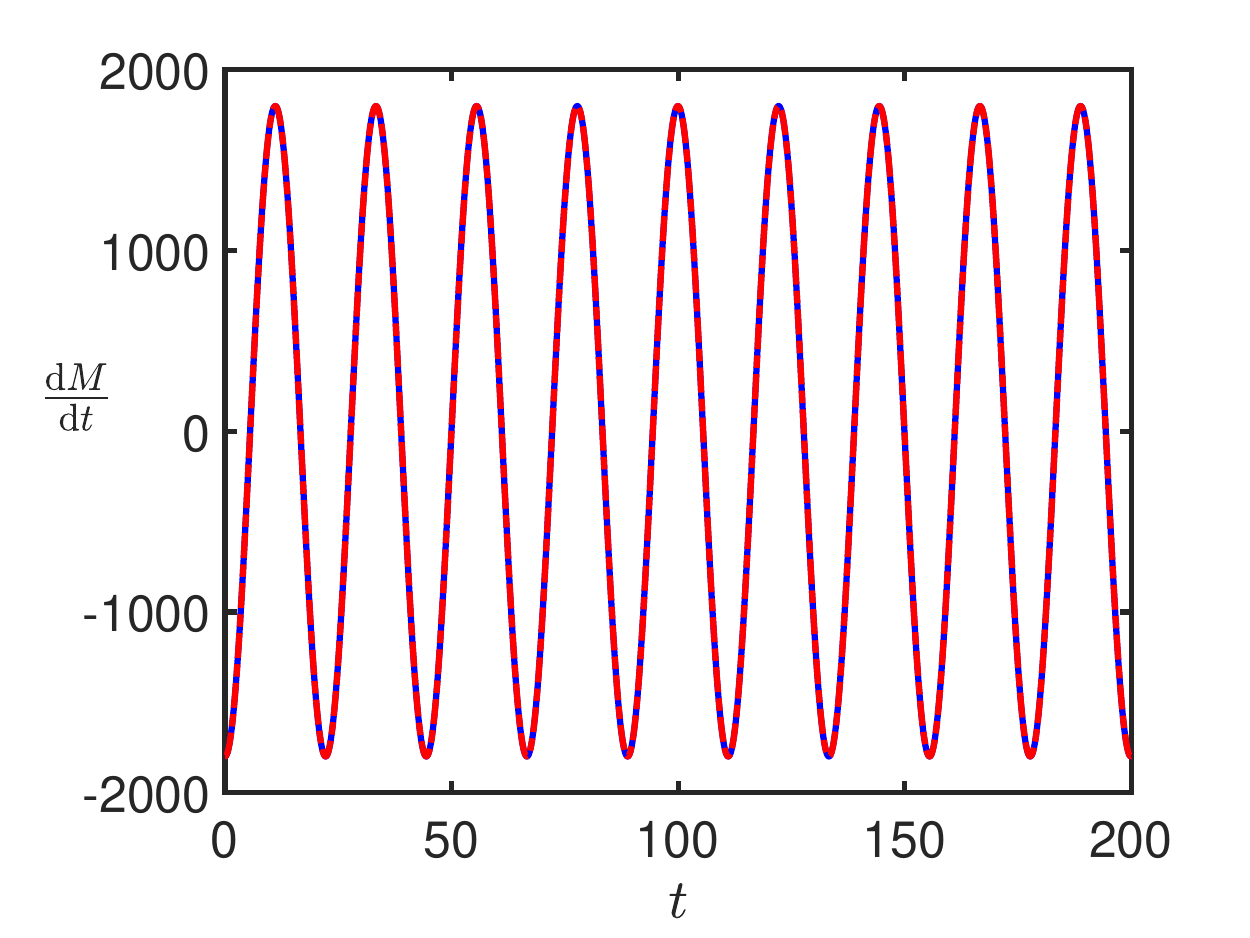}}
	\caption{\small Verification of the generalised conservation laws in the cRB equations, where the blue, solid line is the left-hand side of the relation and the red, dashed line is the right-hand side of the relation, for (a) energy and (b) momentum. Parameters are $L = 300$, $N = 60000$, $c = 1.025$, $\varepsilon = 0.01$, $\alpha = \beta = \delta = \gamma = 1$, $k_1 = k_2 = 1/\sqrt{2}$, $\Delta t = 0.001$ and $\Delta T = \varepsilon \Delta t$.}
	\label{fig:Energy}
	\end{center}
\end{figure}

We now verify the relations \eqref{energy} and \eqref{momentum}. The parameters for the simulation are chosen as $\varepsilon = 0.01$, $c = 1.025$ and $\alpha = \beta = \delta = \gamma = 1$, with solitary wave initial conditions as taken in \eqref{uwIC} with $p = 1$ to provide a significant non-zero mean value. The energy and momentum relations are plotted in Figure \ref{fig:Energy}, where we note that the left-hand side of the mass equation \eqref{mass} is equal to zero to machine precision and therefore is omitted from the plot.

We see that the derivatives of the energy and momentum oscillate with period $f = 2\pi/\omega$ between two different values. To verify the conservation laws, the peaks and troughs were tracked for both the left-hand side and right-hand side of the conservation laws and the absolute percentage error between peaks or troughs was calculated as $4.24 \times 10^{-6}$\% for the energy and $3.60 \times 10^{-10}$\% for the momentum. It is of note that the left-hand side of these laws require a time derivative on discrete data and so the accuracy to which they are conserved is hampered by this limitation. In our case the time step was taken as $\Delta t = 0.001$.

\section{Conclusions}
\label{sec:Conc}
In this paper we developed an asymptotic procedure for the construction of the weakly-nonlinear solution of the Cauchy problem for a system of coupled regularised Boussinesq equations in the high-contrast case when the characteristic speeds in the layers are not close (i.e. the materials of the layers have essentially different elastic properties). Importantly, we removed the need to restrict the initial conditions only to functions with zero mean. The constructed solution was compared to the case when the characteristic speeds in the equations are close (low-contrast case).

We examined the accuracy of the constructed solution numerically, using direct numerical simulations for the coupled Boussinesq equations and our constructed semi-analytical solution, and showed that the constructed solution is a good approximation to the direct numerical simulations, with improving accuracy for the inclusion of additional terms. In particular, we found that the inclusion of terms at $\O{\sqrt{\varepsilon}}$ allows us to take account of the mass with no need for solving additional equations to those at leading order.

We then studied the case of counter-propagating waves within the context of close or distinct characteristic speeds, for solitary wave initial conditions. We showed that radiating solitary waves or wave packets interact almost elastically and their structure appears qualitatively unchanged. Comparison to the case without interaction showed a minor change in amplitude and phase shift for radiating solitary waves, with a more distinct change for the wave packets, suggesting the latter collisions are strongly inelastic.

Next we studied the behaviour of a cnoidal wave initial condition, and solutions for the high-contrast case were compared to that in the low-contrast case. For close characteristic speeds we see that radiating solitary waves are generated from each cnoidal wave peak. The tails interact with the preceding peak, so phase shift and amplitude changes occur in the tail. The amplitude of the main peak is changed by the interactions, however there is no significant phase shift introduced by the interactions. For distinct characteristic speeds, wave packets are generated by each peak, and we can clearly see the evolution of each peak into a wave packet, joined to its neighbours. The individual identities of the wave packets become harder to detect when the wavelength of the cnoidal wave is reduced - the radiation emitted by each wave packet is quickly absorbed by the next one, resulting in a more complicated periodic wave structure.

Finally, we showed that generalised conservation laws can be derived to take account of an initial condition that is not necessarily zero-mean. These conservation laws consist of one local (mass) and two non-local (energy and momentum) relations, with the conserved quantity oscillating with a frequency determined by the evolution of non-zero mass of the initial data, confirmed by numerical simulations. These non-local conservation laws are applicable when the waves propagate on the background of some initially pre-strained basic state characterised by non-zero mass of $u$ and $w$.

The constructed solution can find useful applications in the studies of the scattering of radiating solitary waves by delamination \cite{Khusnutdinova17, Tranter22}, where the previous considerations were restricted to the initial conditions with either small or zero mean value because of the zero-mass contradiction. The propagation of cnoidal waves and undular bores in structures with defects, e.g. \cite{Hooper21, Hooper22}, is an interesting area for future research.

\section*{Acknowledgements}
\label{sec:Ackn}

MRT is grateful to the UK QJMAM Fund for Applied Mathematics for the support of his travel to the ICoNSoM 2019 conference in Rome, Italy, where the initial discussions of this work have taken place. KRK is grateful to the organisers of the Dispersive Hydrodynamics research semester at the Isaac Newton Institute (INI) in the summer of 2022 for the invitation to participate in the program. The work was completed during her stay in Cambridge with the INI support under EPSRC grant number EP/R014604/1.

\begin{appendices}
\section*{Appendices: Numerical Methods}
To solve the equations derived in Section II and its subsections, we make use of pseudo-spectral methods as was done before in \cite{Khusnutdinova19}. We expand on those methods for the case of coupled Boussinesq equations and coupled Ostrovsky equations, as outlined herein.

In the following methods we use the Discrete Fourier Transform (DFT) to calculate the Fourier transform of numerical data. Let us consider a function $u(x,t)$ on a finite domain $x \in [-L, L]$ and we discretise the domain into $N$ equally spaced points, so we have the spacing $\Delta x = 2L/N$. We scale the domain from $x \in [-L, L]$ to $\tilde{x} \in [0, 2\pi]$ via the transform $\tilde{x} = sx + \pi$, where $s = \pi/L$. Denoting $x_{j} = -L + j \Delta x$ for $j=0,\dots,N$, we define the DFT for the function $u(x,t)$ as
\begin{equation*}
	\hat{u} \lb k, t \rb = \frac{1}{\sqrt{N}} \sum_{j=1}^{N} u \lb x_{j}, t \rb e^{-i k x_{j}}, \quad -\frac{N}{2} \leq k \leq \frac{N}{2} - 1,
\end{equation*}
and similarly the IDFT is defined as
\begin{equation*}
	u \lb x, t \rb = \frac{1}{\sqrt{N}} \sum_{k=-N/2}^{N/2 - 1} \hat{u} \lb k, t \rb e^{i k x_{j}}, \quad j =  1, 2, \dots, N,
\end{equation*}
where we have discretised and scaled wavenumber $k \in \mathbb{Z}$. The transforms are implemented using the FFTW3 algorithm \cite{FFTW3}.

% Numerical method for cRB equations
\section{Pseudo-spectral Method for cRB Equations}
\label{sec:cRBNum}
For the coupled Boussinesq equations \eqref{ueq} - \eqref{weq} we use a pseudospectral method similar to the one presented in Ref.~13 where this method was used to solve a single regularised Boussinesq equation. We introduce the change of variables
\begin{equation}
	U = u - \varepsilon u_{xx}, \quad W = w - \varepsilon \beta w_{xx},
	\label{PSUW}
\end{equation}
leading to the modified equations
\begin{align}
	&U_{tt} = u_{xx} + \varepsilon \lsq \frac{1}{2} \lb u^2 \rb_{xx} - \delta \lb u - w \rb \rsq, \notag \\
	&W_{tt} = c^2 w_{xx} + \varepsilon \lsq \frac{\alpha}{2} \lb w^2 \rb_{xx} + \gamma \lb u - w \rb \rsq.
	\label{PScRB}
\end{align}
Note that the variables $U$ and $W$ used here are different from the variables used in the main text.
We take the Fourier transform of \eqref{PSUW}, including the change of domain in $x$, to obtain
\begin{equation}
	\hat{u} = \frac{\hat{U}}{1 + \varepsilon s^2 k^2}, \quad \hat{w} = \frac{\hat{W}}{1 + \varepsilon \beta s^2 k^2}.
	\label{PSuwTransform}
\end{equation}
Similarly, we take the Fourier transform of \eqref{PScRB} and substitute \eqref{PSuwTransform} into this expression to obtain an ordinary differential equation (ODE) in $\hat{u}$ and $\hat{w}$, in the modified domain $\tilde{x}$, taking the form
\begin{align}
	\hat{U}_{tt} &= -\frac{\varepsilon \delta + s^2 k^2}{1 + \varepsilon s^2 k^2} \hat{U} - \frac{\varepsilon s^2 k^2}{2} \mathscr{F} \lset \mathscr{F}^{-1} \lsq \frac{\hat{U}}{1 + \varepsilon s^2 k^2} \rsq^2 \rset + \frac{\varepsilon \delta}{1 + \varepsilon s^2 \beta k^2} \hat{W}, = \hat{S}_1 \lb \hat{U}, \hat{W} \rb, \notag \\
	\hat{W}_{tt} &= -\frac{\varepsilon \gamma + c^2 s^2 k^2}{1 + \varepsilon \beta s^2 k^2} \hat{W} - \frac{\varepsilon \alpha s^2 k^2}{2} \mathscr{F} \lset \mathscr{F}^{-1} \lsq \frac{\hat{W}}{1 + \varepsilon \beta s^2 k^2} \rsq^2 \rset + \frac{\varepsilon \gamma}{1 + \varepsilon s^2 k^2} \hat{U} = \hat{S}_{2} \lb \hat{U}, \hat{W} \rb,
	\label{PSuwODE}
\end{align}
where $\mathscr{F}$ denotes the Fourier transform. We solve this system of ODEs using a 4$^{\mathrm{th}}$-order Runge-Kutta method for time stepping, such as the one used in \cite{Khusnutdinova19, Khusnutdinova19Book}. Therefore we rewrite the system as a series of first-order ODEs, namely
\begin{equation}
	\hat{U}_{t} = \hat{G}, \quad \hat{G}_{t} = \hat{S}_1 \lb \hat{U}, \hat{W} \rb, \quad \hat{W}_{t} = \hat{H}, \quad \hat{H}_{t} = \hat{S}_2 \lb \hat{U}, \hat{W} \rb.
	\label{PScRBSplit}
\end{equation}
We use the discretisation $t = t_n$, $\hat{U}(k, t_{n}) = \hat{U}_n$, $\hat{W}(k, t_{n}) = \hat{W}_n$, $\hat{G}(k, t_{n}) = \hat{G}_{n}$, $\hat{H}(k, t_{n}) = \hat{H}_{n}$ for $n=0,1,2,\dots$, where $t_{n} = n \Delta t$, and $k$ discretises the Fourier space. Taking the Fourier transform of the initial conditions \eqref{uIC} and \eqref{wIC} we obtain initial conditions $\hat{U}_{0}$, $\hat{W}_{0}$, and $\hat{G}_0$, $\hat{H}_{0}$, which are written as
\begin{align}
	\hat{U}_0 &= \lb 1 + \varepsilon s^2 k^2 \rb \mathscr{F} \lset F_{1}(x) \rset, \quad \hat{W}_0 = \lb 1 + \varepsilon \beta s^2 k^2 \rb \mathscr{F} \lset F_{2}(x) \rset, \notag \\
	\hat{G}_0 &= \lb 1 + \varepsilon s^2 k^2 \rb \mathscr{F} \lset V_{1}(x) \rset, \quad \hat{H}_0 = \lb 1 + \varepsilon \beta s^2 k^2 \rb \mathscr{F} \lset V_{2}(x) \rset.
	\label{PSUWIC}
\end{align}
We implement a 4$^{\mathrm{th}}$-order Runge-Kutta method, namely
\begin{align*}
	\hat{U}_{n+1} &= \hat{U}_n + \frac{1}{6} \lsq k_1 + 2 k_2 + 2 k_3 + k_4 \rsq, \quad \hat{G}_{n+1} = \hat{G}_n + \frac{1}{6} \lsq l_1 + 2 l_2 + 2 l_3 + l_4 \rsq, \\
	\hat{W}_{n+1} &= \hat{W}_n + \frac{1}{6} \lsq m_1 + 2 m_2 + 2 m_3 + m_4 \rsq, \quad \hat{H}_{n+1} = \hat{H}_n + \frac{1}{6} \lsq p_1 + 2p_2 + 2p_3 + p_4 \rsq,
\end{align*}
where
\begin{align}
&k_{1} = \Delta t \hat{G}_n, &&l_{1} = \Delta t \hat{S}_{1} \lb \hat{U}_{n}, \hat{W}_{n} \rb, \notag \\
&m_{1} = \Delta t \hat{H}_n, &&p_{1} = \Delta t \hat{S}_{2} \lb \hat{U}_{n}, \hat{W}_{n} \rb, \notag \\
&k_{2} = \Delta t \lb \hat{G}_n + \frac{l_{1}}{2} \rb, &&l_{2} = \Delta t \hat{S}_{1} \lb \hat{U}_{n} + \frac{k_{1}}{2}, \hat{W}_{n} + \frac{m_{1}}{2} \rb, \notag \\
&m_{2} = \Delta t \lb \hat{H}_n + \frac{p_{1}}{2} \rb, &&p_{2} = \Delta t \hat{S}_{2} \lb \hat{U}_{n} + \frac{k_{1}}{2}, \hat{W}_{n} + \frac{m_{1}}{2} \rb, \notag \\
&k_{3} = \Delta t \lb \hat{G}_n + \frac{l_{2}}{2} \rb, &&l_{3} = \Delta t \hat{S}_{1} \lb \hat{U}_{n} + \frac{k_{2}}{2}, \hat{W}_{n} + \frac{m_{2}}{2} \rb, \notag \\
&m_{3} = \Delta t \lb \hat{H}_n + \frac{p_{2}}{2} \rb, &&p_{3} = \Delta t \hat{S}_{2} \lb \hat{U}_{n} + \frac{k_{2}}{2}, \hat{W}_{n} + \frac{m_{2}}{2} \rb, \notag \\
&k_{4} = \Delta t \lb \hat{G}_n + l_{3} \rb, &&l_{4} = \Delta t \hat{S}_{1} \lb \hat{U}_{n} + k_{3}, \hat{W}_{n} + m_{3} \rb, \notag \\
&m_{4} = \Delta t \lb \hat{H}_n + p_{3} \rb, &&p_{4} = \Delta t \hat{S}_{2} \lb \hat{U}_{n} + k_{3}, \hat{W}_{n} + m_{3} \rb.
\label{PSBOstRK4}
\end{align}
To obtain the solution in the real domain we use the relation \eqref{PSuwTransform}. Explicitly we have
\begin{equation}
	u(x,t) = \mathscr{F}^{-1} \lset \frac{\hat{U}}{1 + \varepsilon s^2 k^2} \rset, \quad w(x,t) = \mathscr{F}^{-1} \lset \frac{\hat{W}}{1 + \varepsilon s^2 \beta k^2} \rset.
	\label{PSUWReal}
\end{equation}

% Numerical method for coupled Ostrovsky equations
\section{Pseudo-spectral Method for Ostrovsky equation}
\label{sec:OstNum}
We now consider the solution to the Ostrovsky equation, where we can use a modified Runge-Kutta method as was done in \cite{Khusnutdinova19}. We present an example for the equation governing $\phi_{1}^{-}$, denoted as $\phi$, which can be reduced to the equation for $f^{-}$, denoted $f$. Similarly, we denote $\xi_{-}$ as $\xi$. Explicitly we consider \eqref{phi1eqc1}, as this method can be reduced to solve \eqref{f1Ostc1}, or modified for \eqref{phi2eqc1} and \eqref{f2Ostc1}. We have
\begin{equation}
	\lb 2 \phi_{T} + \lb f \phi \rb_{\xi} + d_{1} \phi_{\xi} + \phi_{\xi \xi \xi} \rb_{\xi} = \delta \phi + H_1(f).
	\label{PSOst}
\end{equation}
We consider the equation on the domains $t \in [0, T]$ and $x \in [-L, L]$. Here $H$ is a function of $f$ and therefore known at the given time step. The form can be found in Section \ref{sec:WNL}. Taking the Fourier transform of \eqref{PSOst} gives (using the transform to $\tilde{\xi}$)
\begin{equation}
	2 \hat{\phi}_{t} + \lb i s k d_1 - i s^3 k^3 \rb \hat{\phi} + i s k \mathscr{F} \lset f \phi \rset = -\frac{i}{sk} \lb \delta \hat{\phi} + \hat{H} \rb.
	\label{PSOstFour}
\end{equation}
Following the method of \cite{Khusnutdinova19} to remove the stiff term from this equation, we multiply through by a multiplicative factor $M$ and introduce a new function $\Phi_{1}$, where $M$ and $\Phi$ are
\begin{equation}
	M = e^{-\frac{i}{2} \lb s^3 k^3 - d_1 s k - \frac{\delta}{sk} \rb t}, \enspace \hat{\Phi} = e^{-\frac{i}{2} \lb s^3 k^3 - d_1 s k - \frac{\gamma}{sk} \rb t} \hat{\phi} = M \hat{\phi}.
	\label{TPhieqOst}
\end{equation}
Substituting this into \eqref{PSOstFour} leads to an ODE for $\Phi$ of the form
\begin{align}
	\hat{\Phi}_{t} = -\frac{i s k}{2} M \mathscr{F} \lset \mathscr{F}^{-1} \lsq \frac{\hat{\Phi}}{M} \rsq \rset - \frac{i}{2sk} M \hat{S}.
	\label{PhiODEOst}
\end{align}
This version of the equation contains less terms than the standard discretisation and therefore we can use an optimised 4$^{\text{th}}$ order Runge-Kutta algorithm. We discretise the time domain as $t_j = j \Delta t$ and the functions as $\hat{\Phi}_{j} = \hat{\Phi} \lb k, t_j \rb$, $\hat{\phi}_{j} = \hat{\phi} \lb k, t_j \rb$ and $\hat{f}_{j} = \hat{f} \lb k, t_j \rb$. Introducing the additional function
\begin{equation}
	E = e^{\frac{i}{4} \lb s^3 k^3 - d_1 sk - \frac{\delta}{sk} \rb \Delta t}, 
	\label{EeqOst}
\end{equation}
we can introduce the optimised Runge-Kutta algorithm in the original variable $\phi$,
	\begin{equation*}
		\hat{\phi}_{j+1} = E_{1}^2 \hat{\phi} + \frac{1}{6} \lsq E_{1}^2 k_1 + 2 E_{1} \lb k_2 + k_3 \rb + k_4 \rsq,
	\end{equation*}
	\vspace{-2em}
	\begin{align}
		k_1 &= -\frac{i s k}{2} \Delta t \mathscr{F} \lset \hat{f}_{j} \mathscr{F}^{-1} \lsq \hat{\phi}_{j} \rsq \rset - \frac{i \Delta t}{2sk} \hat{H}_{1}, &&k_2 = -\frac{i s k}{2} \Delta t \mathscr{F} \lset \hat{f}_{j} \mathscr{F}^{-1} \lsq E \lb \hat{\phi}_{j} + \frac{k_1}{2} \rb \rsq \rset - \frac{i \Delta t}{2sk} \hat{H}, \notag \\
		k_3 &= -\frac{i s k}{2} \Delta t \mathscr{F} \lset \hat{f}_{j} \mathscr{F}^{-1} \lsq E \hat{\phi}_{j} + \frac{k_2}{2} \rsq \rset - \frac{i \Delta t}{2sk} \hat{H}_{1}, &&k_4 = -\frac{i s k}{2} \Delta t \mathscr{F} \lset \hat{f}_{j} \mathscr{F}^{-1} \lsq E^2 \hat{\phi}_{j} + E k_3 \rsq \rset - \frac{i \Delta t}{2sk} \hat{H}.
	\label{RK4OptOst}
	\end{align}
We can apply this algorithm to a homogeneous Ostrovsky equation by setting $\hat{H} = 0$ and replacing the term $f \phi$ with $f^{2}/2$, and therefore it can be applied to all derived Ostrovsky equations in Section \ref{sec:WNL}.

% Numerical method for coupled Ostrovsky equations
\section{Pseudo-spectral Method for coupled Ostrovsky equations}
\label{sec:cOstNum}
We now derive a modified Runge-Kutta method for the coupled Ostrovky equations. Note that these equations are not presented in the main paper, but are in \cite{Khusnutdinova19Book}. Namely, we consider (2.68) and (2.69) for $\phi$, which can be reduced to solve the leading order equations in a similar way to the reduction in Appendix \ref{sec:OstNum}. Explicitly, we will consider
\begin{align}
	\lb 2 \phi_{1 t} + \mu \phi_{1 x} + \lb f_{1} \phi_{1} \rb_{x} + \phi_{1 xxx} \rb_{x} &= \delta \lb \phi_{1} - \phi_{2} \rb + H_{1}, \notag \\
	\lb 2 \phi_{2 t} + \omega \phi_{2 x} + \alpha \lb f_{2} \phi_{2} \rb_{x} + \beta \phi_{2 xxx}  \rb_{x} &= \gamma \lb \phi_{2} - \phi_{1} \rb + H_{2},
	\label{PScOst}
\end{align}
where $\alpha$, $\beta$, $\mu$, $\omega$, $\delta$ and $\gamma$ are constants. Here $H_{1,2}$ are functions of $f_{1,2}$ and, as $f_{1}$, $f_{2}$ are known, this is a known function. For their explicit form, refer to \cite{Khusnutdinova19Book}.  We consider the equation on the domains $t \in [0, T]$ and $x \in [-L, L]$. Taking the Fourier transform of \eqref{PScOst} gives (using the transform to $\tilde{x}$)
	\begin{align}
		2 \hat{\phi}_{1t} + \zeta_1 \hat{\phi}_{1} + i s k \mathscr{F} \lset f_{1} \phi_{1} \rset &= -\frac{i \delta}{sk} \lb \hat{\phi}_{1} - \hat{\phi}_{2} \rb - \frac{i}{sk} \hat{H}_{1}, \notag \\
		2 \hat{\phi}_{2t} + \zeta_2 \hat{\phi}_{2} + i s k \alpha \mathscr{F} \lset f_{2} \phi_{2} \rset &= -\frac{i \gamma}{sk} \lb \hat{\phi}_{2} - \hat{\phi}_{1} \rb - \frac{i}{sk} \hat{H}_{2},
		\label{PScOstFour}
	\end{align}
where 
\begin{equation*}
	\zeta_1 = \lb i s k \mu - i s^3 k^3 \rb, \quad \zeta_2 = \lb i s k \omega - i s^3 k^3 \beta \rb.
\end{equation*}
Following the method of Appendix \ref{sec:OstNum}, to remove the stiff term from this equation we multiply through by a multiplicative factor $M_{1,2}$ and introduce a new function $\Phi_{1,2}$, where $M_{1,2}$ and $\Phi_{1,2}$ take the form
\begin{align}
	M_{1} &= e^{-\frac{1}{2} \lb -\zeta_1 - \frac{i \delta}{sk} \rb t},  &&M_{2} = e^{-\frac{1}{2} \lb -\zeta_2 - \frac{i \gamma}{sk} \rb t}, \notag \\
	\hat{\Phi}_{1} &=  M_{1} \hat{\phi}_{1}, &&\hat{\Phi}_{2} = M_{2} \hat{\phi}_{2}.
	\label{TPhieqcOst}
\end{align}
Substituting this into \eqref{PScOstFour} leads to two ODEs for $\Phi_{1,2}$, which are written as
\begin{align}
	\hat{\Phi}_{1 t} = -\frac{i s k}{2} M_{1} \mathscr{F} \lset \mathscr{F}^{-1} \lsq \frac{\hat{\Phi}_{1}}{M_{1}} \rsq \rset - \frac{i}{2sk} M_{1} \lb \hat{S} - \delta \frac{\Phi_{2}}{M_{2}} \rb, \notag \\
	\hat{\Phi}_{2 t} = -\frac{i \alpha s k}{2} M_{2} \mathscr{F} \lset \mathscr{F}^{-1} \lsq \frac{\hat{\Phi}_{2}}{M_{2}} \rsq \rset - \frac{i}{2sk} M_{2} \lb \hat{S} - \gamma \frac{\Phi_{1}}{M_{1}} \rb.
	\label{PhiODEcOst}
\end{align}
Therefore we can use an optimised 4$^{\text{th}}$ order Runge-Kutta algorithm. We discretise the time domain as $t_j = j \Delta t$ and the functions as $\hat{\Phi}_{i, j} = \hat{\Phi}_{i} \lb k, t_j \rb$, $\hat{\phi}_{i, j} = \hat{\phi}_{i} \lb k, t_j \rb$, $\hat{f}_{i, j} = \hat{f}_{i} \lb k, t_j \rb$, $i = 1,2$.
Introducing the additional functions
\begin{equation}
	E_{1} = e^{\frac{i}{4} \lb s^3 k^3 - \mu sk - \frac{\delta}{sk} \rb \Delta t}, \quad E_{2} = e^{\frac{i}{4} \lb \beta s^3 k^3 - \omega sk - \frac{\gamma}{sk} \rb \Delta t}, 
	\label{EeqcOst}
\end{equation}
we can write the optimised Runge-Kutta algorithm in the original variables $\phi_{1,2}$, taking the form
	\begin{align}
		\hat{\phi}_{1, j+1} &= E_{1}^2 \hat{\phi}_{1} + \frac{1}{6} \lsq E_{1}^2 k_1 + 2 E_{1} \lb k_2 + k_3 \rb + k_4 \rsq, \quad \hat{\phi}_{2, j+1} = E_{2}^2 \hat{\phi}_{1} + \frac{1}{6} \lsq E_{2}^2 l_1 + 2 E_{2} \lb l_2 + l_3 \rb + l_4 \rsq, \notag \\
		k_1 &= -\frac{i s k}{2} \Delta t \mathscr{F} \lset \hat{f}_{1, j} \mathscr{F}^{-1} \lsq \hat{\phi}_{1, j} \rsq \rset - \frac{i \Delta t}{2sk} \lb \hat{H}_{1} - \delta \hat{\phi}_{2, j} \rb, \notag \\
		l_1 &= -\frac{i \alpha s k}{2} \Delta t \mathscr{F} \lset \hat{f}_{2, j} \mathscr{F}^{-1} \lsq \hat{\phi}_{2, j} \rsq \rset - \frac{i \Delta t}{2sk} \lb \hat{H}_{2} - \gamma \hat{\phi}_{1, j} \rb, \notag \\
	k_2 &= -\frac{i s k}{2} \Delta t \mathscr{F} \lset \hat{f}_{1, j} \mathscr{F}^{-1} \lsq E_{1} \lb \hat{\phi}_{1, j} + \frac{k_1}{2} \rb \rsq \rset - \frac{i \Delta t}{2sk} \lb \hat{H}_{1} - \delta E_{2} \lsq \phi_{2, j} + \frac{l_1}{2} \rsq \rb, \notag \\
		l_2 &= -\frac{i \alpha s k}{2} \Delta t \mathscr{F} \lset \hat{f}_{2, j} \mathscr{F}^{-1} \lsq E_{2} \lb \hat{\phi}_{2, j} + \frac{l_1}{2} \rb \rsq \rset - \frac{i \Delta t}{2sk} \lb \hat{H}_{2} - \gamma E_{1} \lsq \phi_{1, j} + \frac{k_1}{2} \rsq \rb, \notag \\
		k_3 &= -\frac{i s k}{2} \Delta t \mathscr{F} \lset \hat{f}_{1, j} \mathscr{F}^{-1} \lsq E_{1} \hat{\phi}_{1, j} + \frac{k_2}{2} \rsq \rset - \frac{i \Delta t}{2sk} \lb \hat{H}_{1} - \delta \lsq E_{2} \hat{\phi}_{2, j} + \frac{l_2}{2} \rsq \rb, \notag \\
		l_3 &= -\frac{i \alpha s k}{2} \Delta t \mathscr{F} \lset \hat{f}_{2, j} \mathscr{F}^{-1} \lsq E_{2} \hat{\phi}_{2, j} + \frac{l_2}{2} \rsq \rset - \frac{i \Delta t}{2sk} \lb \hat{H}_{2} - \gamma \lsq E_{1} \hat{\phi}_{1, j} + \frac{k_2}{2} \rsq \rb, \notag \\
		k_4 &= -\frac{i s k}{2} \Delta t \mathscr{F} \lset \hat{f}_{1, j} \mathscr{F}^{-1} \lsq E_{1}^2 \hat{\phi}_{1, j} + E_{1} k_3 \rsq \rset - \frac{i \Delta t}{2sk} \lb \hat{H}_{1} - \delta \lsq E_{2}^2 \hat{\phi}_{2, j} + E_{2} l_{3} \rsq \rb, \notag \\
		l_4 &= -\frac{i \alpha s k}{2} \Delta t \mathscr{F} \lset \hat{f}_{2, j} \mathscr{F}^{-1} \lsq E_{2}^2 \hat{\phi}_{2, j} + E_{2} l_3 \rsq \rset - \frac{i \Delta t}{2sk} \lb \hat{H}_{2} - \gamma \lsq E_{1}^2 \hat{\phi}_{1, j} + E_{1} k_{3} \rsq \rb.
		\label{RK4OptcOst}
	\end{align}
When calculating the solution using this algorithm, the functions $k_{i}$, $l_{i}$ must be calculated ``in pairs'' as the functions $k_{1}$ and $l_{1}$ are required when evaluating $k_{2}$ and $l_{2}$, and so on. We can apply this algorithm to a homogeneous Ostrovsky equation by setting $\hat{H}_{1,2} = 0$ and replacing the term $f_{i} \phi_{i}$ with $f_{i}^{2}/2$.

\end{appendices}

\bibliographystyle{ieeetr}
\bibliography{NL_KT}

\end{document}